\RequirePackage{fix-cm}
\documentclass[smallextended]{svjour3}       % onecolumn (second format)
\smartqed  % flush right qed marks, e.g. at end of proof
\usepackage{graphicx}
\usepackage{amssymb}
\usepackage[colorlinks=true,linkcolor=blue]{hyperref}
\usepackage{amsmath,amssymb}
\usepackage{graphicx,graphics}% Include figure files
\usepackage{subfigure}
\usepackage{dcolumn}% Align table columns on decimal point
\usepackage{bm}% bold math
\usepackage{float}

\renewcommand{\div}{\mbox{\rm div\,}}
\newcommand{\curl}{\mbox{\rm curl\,}}

\newcommand{\cX}{\mathcal{X}}

\newcommand{\cM}{\mathcal{M}}

\newcommand{\cW}{\mathcal{W}}

\newcommand{\R}{\mathbf{R}}

\begin{document}

\title{Error estimates of deep learning methods for the
nonstationary  Magneto-hydrodynamics equations%\thanks{Grants or other notes
}
%\subtitle{Error analysis for the
%Cahn-Hilliard-Magneto-hydrodynamics problem}

\titlerunning{Error estimates for the  Magneto-hydrodynamics equations}        % if too long for running head

\author{Hailong Qiu
}


\institute{Hailong Qiu \at
School of Mathematics and Physics, Yancheng Institute of Technology,
Yancheng, 224051, China.\\
 This work is supported by the Natural Science Foundation of China (No. 11701498).\\
              \email{qhllf@163.com}           %  \\
}

\date{Received: date / Accepted: date}

\maketitle

\begin{abstract}
In this study, we prove rigourous bounds on the error and stability analysis of deep learning
methods for the nonstationary  Magneto-hydrodynamics equations.
We obtain the approximate ability of the neural network by the convergence of a loss function
and the convergence of a Deep Neural
Network (DNN) to the exact solution.
Moreover, we derive explicit error estimates
 for the solution computed by optimizing the loss function in the DNN approximation of the solution.
\end{abstract}

\section{Introduction}
\label{intro}

Deep learning has been very successfully developed in the artificial
intelligence revolution and forefront of the data science in the last thirty years.
Meanwhile, there has a wide range of applications for deep learning methods in computer vision,
natural language processing and image recognition \cite{Krizhevsky2012,Goodfellow2016,He2016}.
 Recently, as deep neural networks (DNN) are universal function approximators \cite{Hornik1990,Hornik1991,Hornik1997,Xie2011}, it is also natural
to use them for the solutions of partial differential equations (PDEs) \cite{E2018,RaissiP2018,Raissi2018,Sirignano2018,Raissi19,Lu2019,LuK2019,Gulian2019,Yang2020,Cai2022,Dong2022}.
Prominent examples for the application of deep learning methods in PDEs include the deep neural network
approximation of elliptic PDEs \cite{Schwab2019,Kutyniok2021} and nonlinear hyperbolic PDEs \cite{Lye2020,Lye2021}
 and Navier-Stokes equations \cite{Jimenez2018,Wu2018,Fang2020,Lye2020,Mao2020}
 and references therein.

The need of deep learning comes from the fact that
when applying traditional numerical methods in a high dimensional PDEs, the standard methods
may become infeasible. High dimensional PDEs arise in many models
for instance in a variety of contexts such as in derivative
pricing systems, in the financial models,
 credit valuation adjustment problems and  portfolio optimization problems.
These high dimensional nonlinear PDEs are extraordinarily difficult to compute as the
computational effort for traditional approximation methods grows with
the dimension. For example, in  finite
difference methods or finite element methods,
 the number of grids increasing considerably needs
  the memory demands and computational cost
 as the dimension of the
PDEs increases.
However, the deep learning methods in PDEs
present implicit regularization and can surmount the curse of high dimensions \cite{Beck2019,Berner2020}.
In addition, deep learning methods provide a natural framework for bounding unknown
parameters \cite{Raissi2018,Raissi19,Wu2018,Thuerey2020}.

Here our primary goal is on numerical analysis of the deep learning methods for
solving the nonstationary  Magneto-hydrodynamics equations.
 Magneto-hydrodynamics system mainly describes the hydrodynamical behaviors of conducting fluids
subject to external magnetic fields.
The Magneto-hydrodynamics systems are builded by a combination of Navier-Stokes problems and Maxwell problems.
The research of Magneto-hydrodynamics model is a very importance in both
mathematical theory  and practical applications, such as astrophysics,  geophysics,
 the design of cooling systems with liquid metals for a nuclear reactor, meteorology,
   plasma physics and magnetohydrodynamics generators \cite{Moreau1990,2001Davidson}.

So far, there is a lot of literature raising numerical schemes
 applying DNN and machine learning tools for PDEs, including
 the Navier-Stokes equations
%, the Maxwell equations
 and the Magneto-hydrodynamics equations \cite{Nguyen-Thien1999,Raissi2018,Wu2018,Mao2019,Lu2019,Raissi2019,Fang2020,Mao2020,Thuerey2020,Xu2020}.
 The computational algorithm used in deep learning of PDEs \cite{Raissi2018,Gulian2019,Raissi2019}
  involves representing the approximate solution by a DNN,
  in lieu of a finite difference method, spectral method or finite element method,
  and then establishing an appropriate loss function,
   minimizing with such representations,  measuring the
deviation of this representation from the PDEs and the initial and boundary conditions.
It is well known that
optimization of loss functions in a DNN is a non-convex optimization
problem. Therefore, neither the existence nor the uniqueness of a global optimum is
ensured.
The main focus of this paper is considered
 the approximate ability of the neural network by the convergence of a loss function
and the convergence of a DNN to the exact solution.
Furthermore, we obtain explicit error estimates
 for the solution computed by optimizing the loss function in the DNN approximation of the solution. Our method is similar to
the error analysis of machine learning algorithm for the Navier-Stokes equations in \cite{Biswas2022}.

%One important thing to note in this approach is the following. It is well-known that
%optimization of loss functions in a deep neural network is a non-convex optimization
%problem. Therefore, neither the existence nor the uniqueness of a global optimum is
%guaranteed. Typically, repeated application of stochastic gradient descent results in
%reaching a local minimum, which may or may not be global. Nevertheless, we side
%step this issue by taking advantage of the fact that in these applications, the minimum
%value of the loss function at the exact solution must be zero, and obtaining an explicit
%error estimate in terms of the attained value of the loss function. The estimate we
%obtain in turn guarantees that the approximate solution thus constructed converges, in
%the strong topology, to the true solution as the complexity of the networks tends to
%infinity.

The paper is organized as follows.
In section 2, we present some preliminary results,
including the incompressible Magneto-hydrodynamics equations
and neural networks.  In section 3, we obtain
convergence of the loss function and convergence of DNN to the unique solution.
In section 4, we prove some convergence rates and stability of DNN for the
the velocity field and the magnetic field.

\section{\label{Sec2} Preliminaries}

\subsection{\label{Sec2.1} The Magneto-hydrodynamics equations}

 In this paper, we study the deep learning methods
  for the nonstationary magneto-hydrodynamics fluid flow.
  The governing equations are given as follows:
\begin{subequations}\label{eq1.1}
\begin{alignat}{2} \label{eq1.1a} \partial_t \textbf{u}-\nu\Delta\textbf{u}
+(\textbf{u}\cdot \nabla)\textbf{u}\\\nonumber
+S\textbf{B}\times\mbox{curl}\textbf{B}+\nabla p&=\textbf{f},  &&\ \mbox{in} \ \Omega_T:=\Omega\times (0, T],\\
\partial_t \textbf{B}+\mu\mbox{curl}(\mbox{curl} \textbf{B})-\mbox{curl}(\textbf{u}\times \textbf{B})&=0, &&\ \mbox{in}  \ \Omega_T,\\
\mbox{div} \textbf{u}&=0, &&\ \mbox{in} \ \Omega_T,\\
\mbox{div} \textbf{B}&=0, &&\ \mbox{in} \ \Omega_T.
\end{alignat}
\end{subequations}
 The homogeneous boundary conditions and initial conditions are presented:
\begin{subequations}\label{eq2.1}
\begin{alignat}{2} \label{eq2.1a}
 \textbf{u}&=\textbf{0}, \quad &&\ \mbox{on} \ \partial\Omega_T,\\
 \textbf{B}\cdot \textbf{n}&=0, \quad &&\ \mbox{on} \ \partial\Omega_T,\\
\mbox{curl \textbf{B}}\times \textbf{n}&=0, \quad &&\ \mbox{on} \  \partial\Omega_T,
\end{alignat}
\end{subequations}
and
\begin{align}\label{eq3}
\textbf{u}(\textbf{x},0)=\textbf{u}_{0}(\textbf{x}),\ \
\textbf{B}(\textbf{x},0)=\textbf{B}_{0}(\textbf{x}), \ \ \mbox{in}\ \ \Omega,
\end{align}
 where $T>0$ denotes time, and $\Omega \subset \R^2$ is a
bounded and convex domain with continuous boundary $\partial\Omega$.
$\textbf{u}$, $\textbf{B}$ and $p$ denote
the velocity, the magnetic and the pressure, respectively.
 %$\mathbb{D}(\textbf{u}):=\frac{1}{2}(\nabla \textbf{u}+\nabla \textbf{u}^T)$ denotes the strain-rate tensor.
 $\textbf{f}$  is the known body force.
  The positive constants $\nu$ and  $\mu$
stands for the fluid viscous diffusivity coefficient
and the magnetic diffusivity coefficient, respectively.
$S$ denotes the coupling coefficient.

 We introduce some Sobolev spaces
\begin{align*}
&\cX:={H}_0^1(\Omega)^2=\bigl\{\textbf{v}\in H^1(\Omega)^2: \textbf{v}|_{\partial\Omega}=0\bigr\}, \\
& \cW:={H}^1_n(\Omega)^2:=\bigl\{\textbf{w}\in H^1(\Omega)^2: \,\textbf{v}\cdot\textbf{n}|_{\partial\Omega}=0 \bigr\}, \\
& \cM: =L_0^2(\Omega)=\bigl\{ q\in L^2(\Omega), \int_\Omega q d \textbf{x}=0 \bigr\}.
\end{align*}
For convenience, we also define some necessary bilinear terms
\begin{align*}
a_f(\textbf{u},\textbf{v})&=\int_\Omega\nu\nabla\textbf{u}\cdot\nabla\textbf{v}d\textbf{x},\quad
 d(\textbf{v},q)=\int_\Omega q\mbox{div}\textbf{v}d\textbf{x},\\
a_{B}(\textbf{B},\textbf{H})&=\int_\Omega \mu\curl\textbf{B}\cdot\curl\textbf{H}d\textbf{x},
+\int_\Omega \mu\div\textbf{B}\cdot\div\textbf{H}d\textbf{x},
\end{align*}
and trilinear terms
\begin{align*}
b(\textbf{w},\textbf{u},\textbf{v})&= \frac{1}{2}\int_\Omega[(\textbf{w}\cdot\nabla)\textbf{u}]\cdot \textbf{v}
 -[(\textbf{w}\cdot\nabla)\textbf{v}]\cdot \textbf{u} d\textbf{x}=\int_\Omega[(\textbf{w}\cdot\nabla)\textbf{u}]\cdot \textbf{v}
 +\frac{1}{2}[(\nabla\cdot\textbf{w})\textbf{u}]\cdot \textbf{v} d\textbf{x},\\
c_{\widehat{B}}(\textbf{H},\textbf{B},\textbf{v})&=\int_\Omega S\textbf{H}\times \mbox{curl}\textbf{B}\cdot\textbf{v}d\textbf{x},\
\ \quad c_{\widetilde{B}}(\textbf{u},\textbf{B},\textbf{H}) =\int_\Omega (\textbf{u}\times \textbf{B})\cdot\curl\textbf{H}d\textbf{x}.
\end{align*}
Additionally, thanks to integrating by parts, one finds that
\begin{align}\label{eq4}
b(\textbf{u},\textbf{v},\textbf{v})&=0,\  &&\textbf{u}\in \cX, \textbf{v}\in {H}^1(\Omega)^2.
\end{align}
 Employing the identity $(\textbf{B}\times \mbox{curl}\textbf{H},\textbf{v})=(\textbf{v}\times \textbf{H}, \mbox{curl}\textbf{B})$,  it follows that
\begin{align}\label{eq5}
c_{\widehat{B}}(\textbf{B},\textbf{B},\textbf{u})-Sc_{\widetilde{B}}(\textbf{u},\textbf{B},\textbf{B})=0,\ \ \textbf{u}\in \cX,  \textbf{B} \in \cW.
\end{align}
It is noted that the bilinear term $d(\cdot,\cdot)$  satisfies the LBB condition \cite{Girault1986,Temam1983} as follows
\begin{align}\label{eq6}
\sup_{\textbf{v}\in \cX,\textbf{v}\neq \textbf{0}}\dfrac{d(\textbf{v},q)}{\|\textbf{v}\|_{1}}\geq \beta\|q\|,\ \ \  \forall q \in \cM,
\end{align}
where $\beta>0$ is constant depending on $\Omega$.

For simplicity of nations, let us denote
\begin{align*}
\mathfrak{L}_f[\textbf{u},\textbf{B},p]:&=\partial_t \textbf{u}-\nu\Delta\textbf{u}
+(\textbf{u}\cdot \nabla)\textbf{u}\\\nonumber
&\quad+S\textbf{B}\times\mbox{curl}\textbf{B}+\nabla p-\textbf{f}, \\
\mathfrak{L}_B[\textbf{u},\textbf{B}]:&=\partial_t \textbf{B}+\mu\mbox{curl}(\mbox{curl} \textbf{B})-\mbox{curl}(\textbf{u}\times \textbf{B}).
\end{align*}

\subsection{\label{Sec2.2} Neural networks  }

This subsection establishes to approximate the solution of the incompressible Magneto-hydrodynamics equations
with DNN.

Let us denote $\mathfrak{F}_N$ as a DNN with complexity $N$.
Assume $R\in (0,\infty]$, $L\in \mathbb{N}$ and $l_0,\ldots,l_L\in \mathbb{N}$. Let $\sigma:R\rightarrow R$
be a twice differentiable activation function and define
%\begin{align}\label{eq12}
%\Lambda=\Lambda_{L,W,R}:=\bigcup_{\lambda\in\mathbb{N},\lambda\leq L}\bigcup_{l_0,\ldots,l_L\in\{1,\ldots,W\}}\times^{\lambda}_{k=1}\Bigl([-R,R]^{l_k\times l_{k-1}}\times[-R,R]^{l_k} \Bigr).
% \end{align}
the weights and parameters $\theta_k:=(\cW_k,b_k)$ and $A_k:R^{l_{k-1}}\rightarrow R^{l_{k}}:x\mapsto \cW_kx+b_k$ for $1\leq k\leq L$
and we define $f_k^\theta:R^{l_{k-1}}\rightarrow R^{l_{k}}$ by
 \[
f_k^\theta :=\left\{\begin{array}{ll}
 (\sigma\circ A_k^{\theta})(w)  &\, \, \,\,1\leq k<L,\\
  A_k^{\theta}(w) &\,\,  \, \, \,k=L,.
 \end{array}\right.
\]
We define by $(\textbf{u}_\theta,\textbf{B}_\theta):R^{l_0}\times R^{l_0}\rightarrow R^{l_{L}}\times R^{l_{L}}$ the functions that satisfies for all $(w,s)\in  R^{l_{0}}\times R^{l_{0}}$ such that
$$(\textbf{u}_\theta,\textbf{B}_\theta)=\Bigl(\bigl(f^\theta_L\circ f^\theta_{L-1}\circ\cdots \circ f^\theta_{L-1}\bigr)(w),\bigl(f^\theta_L\circ f^\theta_{L-1}\circ\cdots \circ f^\theta_{L-1}\bigr)(s)\Bigl),$$
where in the setting of approximating the Magneto-hydrodynamics equations \eqref{eq1.1}-\eqref{eq3}.

%\begin{theorem} Let $d,n,L,W \in \mathbb{N}$ and let $(\textbf{u}_\theta,\textbf{B}_\theta):R^{d+1}\times R^{d+1}\rightarrow R^{d+1}\times R^{d+1}$ be a neural network solution of \eqref{eq11.1}
%Suppose tthat $\|\sigma\|_{C^n}\geq 1$. Then it holds
%\begin{align}\label{eq12}
%\|\textbf{u}\|_{C^n} +\|\textbf{B}\|_{C^n}\leq 2(16^L)(d+1)^{2n}\bigl(e^2n^4W^3R^n\|\sigma\|_{C^n}\bigr)^{nL}.
% \end{align}
%\end{theorem}

%
\section{\label{Sec3} Convergence of DNN }

The minimization problem of \eqref{eq1.1}-\eqref{eq3} is defined as
 \begin{align}\label{eq7}
 &\inf_{(\textbf{u},\textbf{B},p)\in \mbox{ appropriate Sobolev space}}\Bigl\{\|\mathfrak{L}_f[\textbf{u},\textbf{B},p]\|^2_{L^2(\Omega_T)}
+\|\mathfrak{L}_B[\textbf{u},\textbf{B}]\|^2_{L^2(\Omega_T)}\\\nonumber&\quad
+\|\mbox{div}\textbf{u}\|^2_{L^2(\Omega_T)}+\|\mbox{div}\textbf{B}\|^2_{L^2(\Omega_T)}
+\|\textbf{u}_{|\partial\Omega}\|^2_{L^2(\partial\Omega_T)}+\|\textbf{B}\cdot \textbf{n}_{|\partial\Omega}\|^2_{L^2(\partial\Omega_T)}\Bigr\}.
 \end{align}

In order to approximate $(\textbf{u},\textbf{B},p)$ using a DDN, we consider the following loss function
\begin{align*}
L&=\alpha_1\|\mathfrak{L}_f[\textbf{u}_\theta,\textbf{B}_\theta,p_\theta]\|^2_{L^2(\Omega_T)}+\alpha_2\|\mathfrak{L}_B[\textbf{u}_\theta,\textbf{B}_\theta]\|^2_{L^2(\Omega_T)}
\\\nonumber&\quad+\alpha_3\|\mbox{div}\textbf{u}\|^2_{L^2(\Omega_T)}+\alpha_4\|\mbox{div}\textbf{B}\|^2_{L^2(\Omega_T)}+\alpha_5\|\textbf{u}_\theta|_{\partial\Omega}\|^2_{L^2(\partial\Omega_T)}\\\nonumber&\quad
+\alpha_6\|\textbf{B}_\theta\cdot \textbf{n}|_{\partial\Omega}\|^2_{L^2(\partial\Omega_T)}, \qquad (\textbf{u}_\theta,\textbf{B}_\theta,p_\theta)\in \mathfrak{F}_N.
 \end{align*}

Assume that $\mathfrak{F}_N$ is a finite dimensional function space on a bounded domain $\Omega$. Take a collocation points $\{x_j\}^m_{j=1}\subset \Omega$ and
$\{y_j\}^m_{j=1}\subset \partial\Omega$, find
 \begin{align}\label{eq777}
 &\inf_{(\textbf{u},\textbf{B},p)\in \mathfrak{F}_N}\Bigl\{\alpha_1\sum^m_{j=1}|\mathfrak{L}_f[\textbf{u}(x_j),\textbf{B}(x_j),p(x_j)]|^2
+\alpha_2\sum^m_{j=1}|\mathfrak{L}_B[\textbf{u}(x_j),\textbf{B}(x_j)]|^2\\\nonumber&\quad
+\alpha_3\sum^m_{j=1}|\mbox{div}\textbf{u}(x_j)|^2+\alpha_4\sum^m_{j=1}|\mbox{div}\textbf{B}(x_j)|^2
+\alpha_5\sum^n_{j=1}|\|\textbf{u}(y_j)|_{\partial\Omega}|^2\\\nonumber&\quad
+\alpha_6\sum^n_{j=1}||\textbf{B}(y_j)\cdot \textbf{n}|_{\partial\Omega}|^2\Bigr\}.
 \end{align}
Here, we note that \eqref{eq777} may be use the monte carlo method to compute the corresponding
Lebesgue integrals. Therefore,  let us consider the
 optimization problem as follows:
 \begin{align}\label{eq8}
 &\inf_{(\textbf{u},\textbf{B},p)\in \mathfrak{F}_N}\Bigl\{\|\mathfrak{L}_f[\textbf{u},\textbf{B},p]\|^2_{L^2(\Omega_T)}
+\|\mathfrak{L}_B[\textbf{u},\textbf{B}]\|^2_{L^2(\Omega_T)}
+\|\mbox{div}\textbf{u}\|^2_{L^2(\Omega_T)}\\\nonumber&\quad+\|\mbox{div}\textbf{B}\|^2_{L^2(\Omega_T)}
+\|\textbf{u}|_{\partial\Omega}\|^2_{L^2(\partial\Omega_T)}+\|\textbf{B}\cdot \textbf{n}|_{\partial\Omega}\|^2_{L^2(\partial\Omega_T)}\Bigr\}.
 \end{align}

\subsection{\label{Sec3} Convergence of the loss function }

\begin{lemma} \cite{Biswas2022} Given $\epsilon>0$, assume that $(\textbf{u},\textbf{B},p)$ is the solution of problem \eqref{eq1.1}-\eqref{eq3}, Then there exists $(\textbf{u}_\theta,\textbf{B}_\theta,p_\theta) \in \mathfrak{F}_N$ such that
\begin{align*}
&\sup_{t\in [0,T]}\|\textbf{u}(t)-\textbf{u}_\theta(t)\|_{L^2(\Omega)}\leq C\epsilon,
&&\sup_{t\in [0,T]}\|\textbf{B}(t)-\textbf{B}_\theta(t)\|_{L^2(\Omega)}\leq C\epsilon,\\
&\|\textbf{u}-\textbf{u}_\theta\|_{H^{1,2}(\Omega_T)}\leq C\epsilon,
&&\|\textbf{B}-\textbf{B}_\theta\|_{H^{1,2}(\Omega_T)}\leq C\epsilon,\\
&\|\textbf{u}-\textbf{u}_\theta\|_{L^4([0,T]\times W^{1,4}(\Omega))}\leq C\epsilon,
&&\|\textbf{B}-\textbf{B}_\theta\|_{L^4([0,T]\times W^{1,4}(\Omega))}\leq C\epsilon,\\
&\|p-p_\theta\|_{L^2([0,T]\times H^1(\Omega))}\leq C\epsilon.
 \end{align*}
\end{lemma}

\begin{theorem}  Under the assumptions of \textbf{Lemma 1}, there exists %a  neural network solution
$(\textbf{u}_\theta,\textbf{B}_\theta,p_\theta)\in \mathfrak{F}_N$ such that
  \begin{align*}
 &\inf_{(\textbf{u}_\theta,\textbf{B}_\theta,p_\theta)\in \mathfrak{F}_N}\Bigl\{\|\mathfrak{L}_f[\textbf{u}_\theta,\textbf{B}_\theta,p_\theta]\|^2_{L^2(\Omega_T)}
+\|\mathfrak{L}_B[\textbf{u}_\theta,\textbf{B}_\theta]\|^2_{L^2(\Omega_T)}
+\|\mbox{div}\textbf{u}_\theta\|^2_{L^2(\Omega)}\\\nonumber&\quad+\|\mbox{div}\textbf{B}\|^2_{L^2(\Omega_T)}
+\|\textbf{u}_\theta|_{\partial\Omega}\|^2_{L^2(\partial\Omega_T)}+\|\textbf{B}_\theta\cdot \textbf{n}|_{\partial\Omega}\|^2_{L^2(\partial\Omega_T)}\Bigr\}\leq C\epsilon.
 \end{align*}
\end{theorem}
\noindent \textit{Proof:}\quad Let $(\textbf{u},\textbf{B}, p)$  be the solution of \eqref{eq1.1}-\eqref{eq3}.
%then $\mathfrak{L}_f[\textbf{u},\textbf{B},p]=0$
%and $\mathfrak{L}_B[\textbf{u},\textbf{B}]=0$.
Thus one finds that
\begin{align}\label{eq9}
&\|\mathfrak{L}_f[\textbf{u}_\theta,\textbf{B}_\theta, p_\theta]-\mathfrak{L}_f[\textbf{u},\textbf{B}, p]\|^2_{L^2(\Omega_T)}
+\|\mathfrak{L}_B[\textbf{u}_\theta,\textbf{B}_\theta]-\mathfrak{L}_B[\textbf{u},\textbf{B}]\|^2_{L^2(\Omega_T)}\\\nonumber&\quad
+\|\mbox{div}(\textbf{u}_\theta-\textbf{u})\|^2_{L^2(\Omega_T)}+\|\mbox{div}(\textbf{B}_\theta-\textbf{B})\|^2_{L^2(\Omega_T)}\\\nonumber&\quad
+\|(\textbf{u}_\theta-\textbf{u})|_{\partial\Omega}\|^2_{L^2(\partial\Omega_T)}+\|(\textbf{B}_\theta-\textbf{B})\cdot \textbf{n}|_{\partial\Omega}\|^2_{L^2(\partial\Omega_T)}.
\\\nonumber&=I_1+\ldots+I_6.
 \end{align}
 In the following, we will bound the terms of \eqref{eq9} one by one.
 For nonlinear term, we obtain
  \begin{align}\label{eq76}
 &\quad \|(\textbf{u}\cdot \nabla)\textbf{u}-(\textbf{u}_\theta\cdot \nabla)\textbf{u}_\theta\|^2_{L^2(\Omega_T)}\\\nonumber&
 \leq C\|\textbf{u}-\textbf{u}_\theta\|^2_{L^4(\Omega_T)}\|\textbf{u}_\theta\|^2_{L^4(\Omega_T)}
 \\\nonumber&\quad+C\|\textbf{u}\|^2_{L^4(\Omega_T)}\|\nabla(\textbf{u}-\textbf{u}_\theta)\|^2_{L^4(\Omega_T)},\\\label{eq77}
&\quad \|\textbf{B}\times\mbox{curl}\textbf{B}-\textbf{B}_\theta\times\mbox{curl}\textbf{B}_\theta\|^2_{L^2(\Omega_T)}
\\\nonumber&
 \leq C\|\textbf{B}-\textbf{B}_\theta\|^2_{L^4(\Omega_T)}\|\textbf{B}_\theta\|^2_{L^4(\Omega_T)}
  \\\nonumber&\quad+C\|\textbf{B}\|^2_{L^4(\Omega_T)}\|\nabla(\textbf{B}-\textbf{B}_\theta)\|^2_{L^4(\Omega_T)},\\\label{eq78}
&\quad \|\mbox{curl}(\textbf{u}\times \textbf{B})-\mbox{curl}(\textbf{u}_\theta\times \textbf{B}_\theta)\|^2_{L^2(\Omega_T)}
\\\nonumber&
 \leq C\|\textbf{u}-\textbf{u}_\theta\|^2_{L^4(\Omega_T)}\|\textbf{B}_\theta\|^2_{L^4(\Omega_T)}
  \\\nonumber&\quad+C\|\textbf{u}\|^2_{L^4(\Omega_T)}\|\nabla(\textbf{B}-\textbf{B}_\theta)\|^2_{L^4(\Omega_T)}.
 \end{align}
   Using \textbf{Lemma 1}, it follows that
  \begin{align}\label{eq73}
&\|\partial_t\textbf{u}-\partial_t\textbf{u}_\theta\|^2_{L^2(\Omega_T)}+\|\partial_t\textbf{B}-\partial_t\textbf{B}_\theta\|^2_{L^2(\Omega_T)}\leq C\epsilon^2,\\\label{eq74}
&\|\Delta\textbf{u}-\Delta\textbf{u}_\theta\|^2_{L^2(\Omega_T)}+C\|\curl\curl(\textbf{B}-\textbf{B}_\theta)\|^2_{L^2(\Omega_T)}\leq C\epsilon^2,\\\label{eq75}
&\|\nabla p-\nabla p_\theta\|^2_{L^2(\Omega_T)}\leq C\epsilon^2.
 \end{align}
 Applying the Triangle inequality, we obtain
   \begin{align}\label{eq79}
\|\nabla\textbf{u}_\theta\|^2_{L^4(\Omega_T)}
\leq C\|\nabla\textbf{u}_\theta-\nabla\textbf{u}\|^2_{L^4(\Omega_T)}+\|\nabla\textbf{u}\|^2_{L^4(\Omega_T)},\\\label{eq80}
\|\nabla\textbf{B}_\theta\|^2_{L^4(\Omega_T)}
  \leq C\|\nabla\textbf{B}\|^2_{L^4(\Omega_T)}+\|\nabla(\textbf{B}-\textbf{B}_\theta)\|^2_{L^4(\Omega_T)}.
%\|\nabla\textbf{B}_\theta\|^2_{L^4(\Omega\times[0,T])}
%\leq C\|\nabla(\textbf{B}-\textbf{B}_\theta)\|^2_{L^4(\Omega\times[0,T])}.
 \end{align}
 and
    \begin{align}\label{eq81}
\|\textbf{u}\|^2_{L^4(H^1(\Omega)\times[0,T])}+\|\textbf{B}\|^2_{L^4(H^1(\Omega)\times[0,T])}\leq C.
%\|\nabla\textbf{B}_\theta\|^2_{L^4(\Omega\times[0,T])}
%\leq C\|\nabla(\textbf{B}-\textbf{B}_\theta)\|^2_{L^4(\Omega\times[0,T])}.
 \end{align}
 For term $I_1$, we bound
 \begin{align*}
&\|\mathfrak{L}_f[\textbf{u}_\theta,\textbf{B}_\theta, p_\theta]-\mathfrak{L}_f[\textbf{u},\textbf{B}, p]\|^2_{L^2(\Omega_T)}\\&
\leq \int_0^T\int_{\Omega}|\partial_t(\textbf{u}_\theta-\textbf{u})|^2dxdt+\int_0^T\int_{\Omega}\nu^2|\Delta(\textbf{u}_\theta-\textbf{u})|^2dxdt\\&\quad
+ \int_0^T\int_{\Omega}|(\textbf{u}_\theta\cdot\nabla)\textbf{u}_\theta-(\textbf{u}\cdot\nabla)\textbf{u}|^2dxdt+ \int_0^T\int_{\Omega}|\nabla(p_\theta-p)|^2dxdt\\&\quad
+\int_0^T\int_{\Omega}S^2|\textbf{B}_\theta\times\curl\textbf{B}_\theta-\textbf{B}\times\curl\textbf{B}|^2dxdt\\&
\leq C\epsilon.
 \end{align*}
For term $I_2$, we estimate
  \begin{align*}
&\quad\|\mathfrak{L}_B[\textbf{u}_\theta,\textbf{B}_\theta]-\mathfrak{L}_B[\textbf{u},\textbf{B}]\|^2_{L^2(\Omega_T)}\\&
\leq  \int_0^T\int_{\Omega}|\partial_t(\textbf{B}_\theta-\textbf{B})|^2dxdt+\int_0^T\int_{\Omega}\mu^2|\curl\curl(\textbf{B}_\theta-\textbf{B})|^2dxdt\\&\quad
+ \int_0^T\int_{\Omega}|\curl(\textbf{u}_\theta\times\textbf{B}_\theta)-\curl(\textbf{u}\times\textbf{B})|^2dxdt\\&
\leq C\epsilon.
 \end{align*}
 It is similar to $I_3+I_4+I_5+I_6$, we have
 \begin{align*}
I_3+I_4+I_5+I_6\leq C\epsilon.
 \end{align*}
 The desired result is derived.
 The proof is completed.
 $$\eqno\Box$$

\subsection{\label{Sec3} Convergence of DNN to the unique solution }

 \begin{theorem} Assume that $(\textbf{u},\textbf{B}, p) \in \cX\times \cW\times \cM$  is a unique solution to \eqref{eq1.1}-\eqref{eq3}, Then
when the sequence $(\textbf{u}^n_\theta,\textbf{B}^n_\theta, p^n_\theta)$ of problem \eqref{eq8} is uniformly bounded and equicontinuous,
 the neural networks solution $(\textbf{u}^n_\theta,\textbf{B}^n_\theta, p^n_\theta)$  converges strongly to $(\textbf{u},\textbf{B}, p)$.
 % Furthermore, $(\textbf{u}^n_\theta,\textbf{B}^n_\theta, p^n_\theta)$ uniformly converges to $(\textbf{u},\textbf{B}, p)$ in $\Omega$, respectively.
\end{theorem}
\noindent \textit{Proof:}\quad Let $(\textbf{u}^n_\theta,\textbf{B}^n_\theta, p^n_\theta) \in \mathfrak{F}_N$ is the solution of problem \eqref{eq8}, then $(\textbf{u}^n_\theta,\textbf{B}^n_\theta, p^n_\theta)$
satisfy
\begin{subequations}\label{eq10.1}
\begin{alignat}{2} \label{eq10.1a} &(\frac{d\textbf{u}^n_\theta}{dt},\textbf{v})+a(\textbf{u}^n_\theta,\textbf{v})
+b(\textbf{u}^n_\theta,\textbf{u}_\theta,\textbf{v})\\\nonumber&\quad
+Sc_{\widehat{B}}(\textbf{B}^n_\theta,\textbf{B}^n_\theta,\textbf{v})-d(\textbf{v},p_\theta)+d(\textbf{u}^n_\theta,q)=(\textbf{f},\textbf{v}), \\
&(\frac{d\textbf{B}^n_\theta}{d t},\textbf{H})+a_B(\textbf{B}^n_\theta,\textbf{H})
-c_{\widetilde{B}}(\textbf{u}^n_\theta,\textbf{B}^n_\theta,\textbf{H})=0,
\end{alignat}
\end{subequations}
for all $(\textbf{v},\textbf{H},q)\in \mathfrak{F}_N$.

Taking $(\textbf{v},\textbf{H},q)=(\textbf{u}^n_\theta,\textbf{B}^n_\theta, p^n_\theta)$ in \eqref{eq10.1},  it follows that
 \begin{align*}
\frac{d}{dt}\bigl(\|\textbf{u}^n_\theta\|^2_{L^2(\Omega)}+S\|\textbf{B}^n_\theta\|^2_{L^2(\Omega)}\bigr)+\nu\|\nabla\textbf{u}^n_\theta\|^2_{L^2(\Omega)}
+c_0\mu S\|\nabla\textbf{B}^n_\theta\|^2_{L^2(\Omega)}\leq C\|\textbf{f}\|^2_{L^2(\Omega)}.
 \end{align*}
 Using the Gr\"{o}nwall inequality, we can get
 \begin{align}\label{eq11}
\sup_{0\leq t\leq T}\bigl(\|\textbf{u}^n_\theta(t)\|^2_{L^2(\Omega)}&+\|\textbf{B}^n_\theta(t)\|^2_{L^2(\Omega)}\bigr)+\int^T_{0}\|\nabla\textbf{u}^n_\theta\|^2_{L^2(\Omega)}dt
\\&\nonumber+\int^T_{0}\|\nabla\textbf{B}^n_\theta\|^2_{L^2(\Omega)}dt\leq C.
 \end{align}

Setting $(\textbf{v},\textbf{H}, q)=(-\Delta\textbf{u}^n_\theta,\curl\curl\textbf{B}^n_\theta, 0)$ in \eqref{eq10.1}, we obtain
\begin{align}\label{eq12}
&\frac{1}{2}\frac{d}{dt}\bigl(\|\nabla\textbf{u}^n_\theta\|^2_{L^2(\Omega)}+Sc_0\|\nabla\textbf{B}^n_\theta\|^2_{L^2(\Omega)}\bigr)
+\nu\|A\textbf{u}^n_\theta\|^2_{L^2(\Omega)}
+\mu S\|\Delta\textbf{B}^n_\theta\|^2_{L^2(\Omega)}\\\nonumber&\quad
\leq |(\textbf{f},A\textbf{u}^n_\theta)
+|b(\textbf{u}^n_\theta,\textbf{u}^n_\theta,A\textbf{u}_\theta)|
+|Sc_{\widehat{B}}(\textbf{B}^n_\theta,\textbf{B}^n_\theta,A\textbf{u}^n_\theta)|
+|Sc_{\widetilde{B}}(\textbf{u}^n_\theta,\textbf{B}^n_\theta,-\Delta\textbf{B}^n_\theta)|.
 \end{align}
 By the H\"{o}lder, Young and Sobolev inequalities, we have
 \begin{align*}
 |(\textbf{f},A\textbf{u}^n_\theta)|&\leq \frac{\nu}{8}\|A\textbf{u}^n_\theta\|^2_{L^2(\Omega)}+C\|\textbf{f}\|^2_{L^2(\Omega)},\\
|b(\textbf{u}_\theta,\textbf{u}^n_\theta,\textbf{u}_\theta)|&\leq \|\textbf{u}^n_\theta\|^{\frac{1}{2}}_{L^2(\Omega)}\|\nabla\textbf{u}^n_\theta\|_{L^2(\Omega)}\|A\textbf{u}^n_\theta\|^{\frac{3}{2}}_{L^2(\Omega)}\\&
\leq \frac{\nu}{8}\|A\textbf{u}^n_\theta\|^2_{L^2(\Omega)}+C\|\textbf{u}^n_\theta\|^2_{L^2(\Omega)}\|\nabla\textbf{u}^n_\theta\|^4_{L^2(\Omega)},\\
|Sc_{\widehat{B}}(\textbf{B}^n_\theta,\textbf{B}_\theta,A\textbf{u}^n_\theta)|&\leq C\|\textbf{B}^n_\theta\|_{L^4(\Omega)}\|\curl\textbf{B}^n_\theta\|_{L^4(\Omega)}\|A\textbf{u}^n_\theta\|_{L^2(\Omega)}\\&
\leq C\|\textbf{B}^n_\theta\|^{\frac{1}{2}}_{L^2(\Omega)}\|\textbf{B}^n_\theta\|^{\frac{1}{2}}_{H^1(\Omega)}\|\curl\textbf{B}^n_\theta\|_{L^4(\Omega)}\|A\textbf{u}^n_\theta\|_{L^2(\Omega)}\\&
\leq \frac{\nu}{8}\|A\textbf{u}^n_\theta\|^2_{L^2(\Omega)}+\frac{\mu S}{8}\|\Delta\textbf{B}^n_\theta\|^2_{L^2(\Omega)}\\&\quad
+\|\textbf{B}^n_\theta\|^2_{L^2(\Omega)}\|\textbf{B}^n_\theta\|^2_{H^1(\Omega)}\|\nabla\textbf{B}^n_\theta\|_{L^2(\Omega)},\\
|Sc_{\widetilde{B}}(\textbf{u}^n_\theta,\textbf{B}^n_\theta,-\Delta\textbf{B}^n_\theta)|&\leq C\bigl(\|\textbf{u}^n_\theta\|_{L^4(\Omega)}\|\nabla\textbf{B}^n_\theta\|_{L^4(\Omega)}
+\|\nabla\textbf{u}^n_\theta\|_{L^4(\Omega)}\|\textbf{B}^n_\theta\|_{L^4(\Omega)}\bigr)\|\Delta\textbf{B}^n_\theta\|_{L^2(\Omega)}\\&
\leq \frac{\nu}{8}\|A\textbf{u}^n_\theta\|^2_{L^2(\Omega)}+\frac{\mu S}{8}\|\Delta\textbf{B}^n_\theta\|^2_{L^2(\Omega)}
+C\bigl(\|\textbf{u}^n_\theta\|^2_{L^2(\Omega)}\|\nabla\textbf{u}_\theta\|^2_{L^2(\Omega)}\\&\quad
+\|\textbf{B}^n_\theta\|^2_{L^2(\Omega)}\|\nabla\textbf{B}^n_\theta\|^2_{L^2(\Omega)}\bigr)\bigl(\|\nabla\textbf{u}^n_\theta\|_{L^2(\Omega)}+\|\nabla\textbf{B}^n_\theta\|_{L^2(\Omega)}\bigr).
 \end{align*}
 Combining the above inequalities with \eqref{eq12}, we obtain
 \begin{align*}
&\frac{d}{dt}\bigl(\|\nabla\textbf{u}^n_\theta\|^2_{L^2(\Omega)}+Sc_0\|\nabla\textbf{B}^n_\theta\|^2_{L^2(\Omega)}\bigr)
+\nu\|A\textbf{u}^n_\theta\|^2_{L^2(\Omega)}
+\mu S\|\Delta\textbf{B}^n_\theta\|^2_{L^2(\Omega)}\\\nonumber&
\leq C\|\textbf{f}\|^2_{L^2(\Omega)}+\bigl(\|\textbf{u}^n_\theta\|^2_{L^2(\Omega)}\|\nabla\textbf{u}^n_\theta\|^2_{L^2(\Omega)}
+\|\textbf{B}^n_\theta\|^2_{L^2(\Omega)}\|\nabla\textbf{B}^n_\theta\|^2_{L^2(\Omega)}\bigr)\\\nonumber&\quad
\times\bigl(\|\nabla\textbf{u}^n_\theta\|^2_{L^2(\Omega)}
+\|\nabla\textbf{B}^n_\theta\|^2_{L^2(\Omega)}\bigr).
 \end{align*}
 With employing the Gr\"{o}nwall inequality, it follows that
  \begin{align}\label{eq13}
\sup_{0\leq t\leq T}\bigl(\|\textbf{u}^n_\theta(t)\|^2_{H^1(\Omega)}&+\|\textbf{B}^n_\theta(t)\|^2_{H^1(\Omega)}\bigr)+\int^T_{0}\|\textbf{u}^n_\theta\|^2_{H^2(\Omega)}dt
\\&\nonumber+\int^T_{0}\|\textbf{B}^n_\theta\|^2_{H^2(\Omega)}dt\leq C.
 \end{align}

 Differentiate both sides of \eqref{eq10.1} with respect to $t$, we have
 \begin{subequations}\label{eq14.1}
\begin{alignat}{2} \label{eq14.1a} &(\frac{d^2\textbf{u}^n_\theta}{dt^2},\textbf{v})+a(\frac{d\textbf{u}^n_\theta}{dt},\textbf{v})
+b(\frac{d\textbf{u}^n_\theta}{dt},\textbf{u}_\theta,\textbf{v})+b(\textbf{u}^n_\theta,\frac{d\textbf{u}^n_\theta}{dt},\textbf{v})\\\nonumber&\quad
+Sc_{\widehat{B}}(\frac{d\textbf{B}^n_\theta}{dt},\textbf{B}_\theta,\textbf{v})+Sc_{\widehat{B}}(\textbf{B}^n_\theta,\frac{d\textbf{B}^n_\theta}{dt},\textbf{v})
-d(\textbf{v},\frac{dp^n_\theta}{dt})%\\\nonumber&\quad
+d(\frac{d\textbf{u}^n_\theta}{dt},q)=(\textbf{f}_t,\textbf{v}), \\
&(\frac{d^2\textbf{B}^n_\theta}{dt^2},\textbf{H})+a_B(\frac{d\textbf{B}^n_\theta}{dt},\textbf{H})
-c_{\widetilde{B}}(\frac{d\textbf{u}^n_\theta}{dt},\textbf{B}^n_\theta,\textbf{H})-c_{\widetilde{B}}(\textbf{u}^n_\theta,\frac{d\textbf{B}^n_\theta}{dt},\textbf{H})=0,
\end{alignat}
\end{subequations} for all $(\textbf{v},\textbf{H},q)\in \mathfrak{F}_N$.

Taking $(\textbf{v},\textbf{H},q)=(\frac{d\textbf{u}^n_\theta}{dt},\frac{d\textbf{B}^n_\theta}{dt}, \frac{dp^n_\theta}{dt})$ in \eqref{eq14.1},  one finds that
\begin{align}\label{eq15}
&\frac{1}{2}\frac{d}{dt}\Bigl(\Bigl\|\frac{d\textbf{u}^n_\theta}{dt}\Bigl\|^2_{L^2(\Omega)}+S\Bigl\|\frac{d\textbf{B}^n_\theta}{dt}\Bigr\|^2_{L^2(\Omega)}\Bigr)+\nu\Bigl\|\nabla\frac{d\textbf{u}^n_\theta}{dt}\Bigr\|^2_{L^2(\Omega)}
+\mu S c_0\Bigl\|\nabla\frac{d\textbf{B}^n_\theta}{d t}\Bigr\|^2_{L^2(\Omega)}\\\nonumber&\quad
\leq |(\textbf{f}_t, \frac{d\textbf{u}^n_\theta}{d t})|
+|b(\frac{d\textbf{u}^n_\theta}{d t},\textbf{u}^n_\theta,\frac{d\textbf{u}_\theta}{d t})|
+|Sc_{\widehat{B}}(\textbf{B}^n_\theta,\frac{d\textbf{B}^n_\theta}{d t},\frac{d\textbf{u}^n_\theta}{d t})|
+|Sc_{\widetilde{B}}(\textbf{u}^n_\theta,\frac{d\textbf{B}^n_\theta}{d t},\frac{d\textbf{B}^n_\theta}{d t})|.
 \end{align}
 Making using of the H\"{o}lder, Young and Embedded inequalities, we derive
 \begin{align*}
|(\textbf{f}_t, \frac{d\textbf{u}^n_\theta}{d t})|&\leq C\|\textbf{f}_t\|^2_{L^2(\Omega)}+\frac{\nu}{6}\Bigl\|\frac{d\textbf{u}^n_\theta}{d t}\Bigl\|^2_{L^2(\Omega)},\\
|b(\frac{d\textbf{u}^n_\theta}{d t},\textbf{u}_\theta,\frac{d\textbf{u}_\theta}{d t})|&\leq C\Bigl\|\frac{d\textbf{u}^n_\theta}{d t}\Bigl\|_{L^4(\Omega)}\Bigl\|\frac{d\textbf{u}^n_\theta}{d t}\Bigl\|_{L^2(\Omega)}\Bigl\| \nabla\textbf{u}^n_\theta\|_{L^4(\Omega)}\\&
\leq \frac{\nu}{6}\Bigl\|\nabla\frac{d\textbf{u}^n_\theta}{d t}\Bigl\|^2_{L^2(\Omega)}+C\Bigl\|\frac{d\textbf{u}^n_\theta}{d t}\Bigl\|_{L^2(\Omega)}\Bigl\|\textbf{u}_\theta\|_{H^2(\Omega)},\\
|Sc_{\widehat{B}}(\textbf{B}_\theta,\frac{d\textbf{B}_\theta}{d t},\frac{d\textbf{u}_\theta}{d t})|&\leq C\Bigl\|\frac{d\textbf{u}_\theta}{d t}\Bigl\|_{L^4(\Omega)}\Bigl\|\frac{d\textbf{B}_\theta}{d t}\Bigl\|_{L^2(\Omega)}\Bigl\| \nabla\textbf{B}^n_\theta\|_{L^4(\Omega)}\\&
 \leq \frac{\nu}{6}\Bigl\|\nabla\frac{d\textbf{u}^n_\theta}{d t}\Bigl\|^2_{L^2(\Omega)}+C\Bigl\|\frac{d\textbf{B}^n_\theta}{d t}\Bigl\|_{L^2(\Omega)}\Bigl\|\textbf{B}^n_\theta\|_{H^2(\Omega)},\\
|Sc_{\widetilde{B}}(\textbf{u}^n_\theta,\frac{d\textbf{B}^n_\theta}{d t},\frac{d\textbf{B}^n_\theta}{d t})|&\leq C\|\textbf{u}^n_\theta\|_{L^\infty(\Omega)}\Bigl\|\frac{d\textbf{B}^n_\theta}{d t}\Bigl\|_{L^2(\Omega)}\Bigl\| \nabla\frac{d\textbf{B}^n_\theta}{d t}\|_{L^2(\Omega)}\\&
 \leq \frac{\mu S}{6}\Bigl\|\nabla\frac{d\textbf{B}^n_\theta}{d t}\Bigl\|^2_{L^2(\Omega)}+C\Bigl\|\frac{d\textbf{B}^n_\theta}{d t}\Bigl\|_{L^2(\Omega)}\Bigl\|\textbf{u}^n_\theta\|_{H^2(\Omega)}.
 \end{align*}
  Combining the above inequalities with \eqref{eq15} yields
  \begin{align*}
&\frac{d}{dt}\Bigl(\Bigl\|\frac{d\textbf{u}^n_\theta}{dt}\Bigl\|^2_{L^2(\Omega)}+S\Bigl\|\frac{d\textbf{B}^n_\theta}{dt}\Bigr\|^2_{L^2(\Omega)}\Bigr)+\nu\Bigl\|\nabla\frac{d\textbf{u}^n_\theta}{dt}\Bigr\|^2_{L^2(\Omega)}
+\mu S c_0\Bigl\|\nabla\frac{d\textbf{B}^n_\theta}{d t}\Bigr\|^2_{L^2(\Omega)}\\\nonumber&\quad
\leq  C\|\textbf{f}_t\|^2_{L^2(\Omega)}+C\bigl(\|\textbf{u}^n_\theta\|^2_{H^2(\Omega)}+\|\textbf{B}^n_\theta\|^2_{H^2(\Omega)}\bigr)\Bigl(\Bigl\|\frac{d\textbf{u}^n_\theta}{dt}\Bigl\|^2_{L^2(\Omega)}
+S\Bigl\|\frac{d\textbf{B}^n_\theta}{dt}\Bigr\|^2_{L^2(\Omega)}\Bigr).
 \end{align*}
By applying the Gr\"{o}nwall inequality, we obtain
  \begin{align}\label{eq16}
\sup_{0\leq t\leq T}\Bigl(\Bigl\|\frac{d\textbf{u}^n_\theta}{dt}\Bigl\|^2_{L^2(\Omega)}&+\Bigl\|\frac{d\textbf{B}^n_\theta}{dt}\Bigr\|^2_{L^2(\Omega)}\Bigr)
+\int^T_{0}\Bigl\|\nabla\frac{d\textbf{u}^n_\theta}{dt}\Bigr\|^2_{L^2(\Omega)}dt\\&\nonumber
+\int^T_{0}\Bigl\|\nabla\frac{d\textbf{B}^n_\theta}{d t}\Bigr\|^2_{L^2(\Omega)}dt\leq C.
 \end{align}
 From  \eqref{eq11}, \eqref{eq13} and  \eqref{eq16},  we can obtain that $\{\textbf{u}^n_\theta\},\{\frac{d\textbf{u}^n_\theta}{dt}\},
 \{\textbf{B}^n_\theta\}$ and $\{\frac{d\textbf{B}^n_\theta}{dt}\}$ are uniformly bounded in
 $L^2([0,T], H^1(\Omega))$.
Applying the Aubin-Lions's compactness lemma, there exists a subsequence of
$\{\textbf{u}^n_\theta\},\{\frac{d\textbf{u}^n_\theta}{dt}\},
 \{\textbf{B}^n_\theta\}$ and $\{\frac{d\textbf{B}^n_\theta}{dt}\}$ (still
 denoted by the $\{\textbf{u}^n_\theta\},\{\frac{d\textbf{u}^n_\theta}{dt}\},
 \{\textbf{B}^n_\theta\}$ and $\{\frac{d\textbf{B}^n_\theta}{dt}\}$), which converges to
$\textbf{u}\in L^{\infty}([0,T], L^2(\Omega))\cap L^{2}([0,T], \cX)$ and
$\textbf{B}\in L^{\infty}([0,T], L^2(\Omega))\cap L^{2}([0,T], \cW)$ such that
$$\textbf{u}^n_\theta\rightarrow \textbf{u} \qquad \mbox{in} \qquad L^{2}([0,T], L^2(\Omega)),$$
and
$$\textbf{B}^n_\theta\rightarrow \textbf{B}\qquad \mbox{in} \qquad L^{2}([0,T], L^2(\Omega)).$$

Finally, we have been ready for passing to the limit as $n\rightarrow\infty$ in the weak sense, it
is not difficult to show that $(\textbf{u},\textbf{B})$ satisfy \eqref{eq10.1} in a weak formulation. This result presents that
$\textbf{u}^n_\theta$ and $\textbf{B}^n_\theta$ converge strongly to $\textbf{u}$ and $\textbf{B}$ in $L^{2}([0,T], L^2(\Omega))$.
Similarly, using \eqref{eq6}, we can derive $p^n_\theta$ converge strongly to $p$ in $L^{2}([0,T], L^2(\Omega))$.
 The desired result is derived.
 The proof is completed.
 $$\eqno\Box$$

%%%%%%%%%%%%%%%%%%%%%%%%%%%%%%%%%%%%%%%%%%%%%%%%%%%%%%%%%%%%%%%%%%%%%%%%%%%%%%%%%%%%%%%%%%%%%%%%%%%%%%%%%%%%%%%%%%%%

\section{\label{Sec4} Convergence rates and stability of DNN }

In this section, we drive some convergence rates and stability of DNN  for problem \eqref{eq1.1}-\eqref{eq3}.
 Assume that a DNN solution $(\textbf{u}_\theta,\textbf{B}_\theta, p_\theta)$ of problem \eqref{eq8} satisfy
 \begin{align}\label{eq17}
\frac{d\textbf{u}_\theta}{dt}+\nu A_f\textbf{u}_\theta+B[\textbf{u}_\theta,\textbf{u}_\theta]+SC_f[\textbf{B}_\theta,\textbf{B}_\theta]=\mathbb{P}\textbf{f},\\
\frac{d\textbf{B}_\theta}{dt}+\mu A_B\textbf{B}_\theta-C_B[\textbf{u}_\theta,\textbf{B}_\theta]=0,
 \end{align}
 where $\mathbb{P}$ and $\mathbb{Q}$ are Leray projections, $A_f:=-\mathbb{P}\Delta$ and $A_B:=\mathbb{Q}\curl\curl$ are Stokes operator and Maxwell operator, respectively.
 $B[\textbf{u},\textbf{u}]$, $C_f[\textbf{B},\textbf{B}]$ and $C_B[\textbf{u},\textbf{B}]$ are defined as follows:
  \begin{align*}
 B[\textbf{u},\textbf{u}]:=\mathbb{P}[(\textbf{u}\cdot\nabla)\textbf{u}],\quad C_f[\textbf{B},\textbf{B}]:=\mathbb{P}[\curl(\textbf{B}\times\textbf{B}],\quad C_B[\textbf{u},\textbf{B}]:=\mathbb{Q}[\curl(\textbf{u}_\theta\times\textbf{B}_\theta)].
 \end{align*}

Here we use the similar  technique \cite{Biswas2022}, thus we introduce the Hodge decomposition.
The main idea of Hodge decomposition is to decompose a vector $\textbf{w}\in L^2(\Omega)$ uniquely
into a divergence-free part $\textbf{w}^1$ and an irrotational part  $\textbf{w}^2$, which is orthogonal in $L^2(\Omega)$
to $\textbf{w}^1$, i.e.,
  \begin{align}\label{eq333}
\textbf{w}=\textbf{w}^1+\textbf{w}^2,\qquad \nabla\cdot \textbf{w}^1=0, \qquad (\textbf{w}^1,\textbf{w}^2)=0.
 \end{align}
 Consider an approximate solution $(\textbf{u}_\theta,\textbf{B}_\theta, p_\theta) \in \mathfrak{F}_N$ and denote
 \begin{align}\label{eq38}
\mathfrak{L}_f[\textbf{u}_\theta,\textbf{B}_\theta, p_\theta]&=\widehat{\textbf{f}},\\\label{eq39}
\mathfrak{L}_B[\textbf{u}_\theta,\textbf{B}_\theta]&=\textbf{0},\\\label{eq40}
\nabla\cdot\textbf{u}_\theta&=g,\\\label{eq41}
\nabla\cdot\textbf{B}_\theta&=h.
 \end{align}
 Applying the Hodge decomposition on $(\textbf{u}_\theta,\textbf{B}_\theta)$, we have
 \begin{align}\label{eq42}
\textbf{u}_\theta&=\mathbb{P}\textbf{u}_\theta+(\mathbb{I}-\mathbb{P})\textbf{u}_\theta=:\textbf{u}^1_\theta+\textbf{u}^2_\theta,\\\label{eq43}
\textbf{B}_\theta&=\mathbb{Q}\textbf{B}_\theta+(\mathbb{I}- \mathbb{Q})\textbf{B}_\theta=:\textbf{B}^1_\theta+\textbf{B}^2_\theta,
 \end{align}

\begin{theorem} Assume that $(\textbf{u},\textbf{B},p)$ is a strong solution of problem \eqref{eq1.1}-\eqref{eq3}
 and $(\textbf{u}_\theta,\textbf{B}_\theta,p_\theta)$ such that
 \begin{align}\label{eq18}
\frac{d\textbf{u}^1_\theta}{dt}+\nu A_f\textbf{u}^1_\theta+B[\textbf{u}^1_\theta,\textbf{u}^1_\theta]+SC_f[\textbf{B}^1_\theta,\textbf{B}^1_\theta]=\mathbb{P}\textbf{f}+\Lambda,\\\label{eq19}
\frac{d\textbf{B}^1_\theta}{dt}+\mu A_B\textbf{B}^1_\theta-C_B[\textbf{u}^1_\theta,\textbf{B}^1_\theta]=\Pi,
 \end{align}
 where
  \begin{align}\label{eq20}
\int_{0}^T\|\Lambda\|^{2}_{H^{-1}}dt+\int_{0}^T\|\Pi\|^{2}_{\cW^{-1}}dt\leq C\epsilon,
 \end{align}
 and assume
  \begin{align}\label{eq21}
\|\textbf{u}^1_{\theta,0}-\textbf{u}_{0}\|^{2}_{L^2(\Omega)}+\|\textbf{B}^1_{\theta,0}-\textbf{B}_{0}\|^{2}_{L^2(\Omega)}\leq C\epsilon,
 \end{align}
then the following bound satisfy
 \begin{align}\label{eq22}
\sup_{t \in [0,T]}\bigl(\|\textbf{u}(t)-\textbf{u}^1_{\theta}(t)\|^2_{L^2(\Omega)}+\|\textbf{B}(t)-\textbf{B}^1_{\theta}(t)\|^2_{L^2(\Omega)}\bigr)\leq C\epsilon.
 \end{align}
\end{theorem}
\noindent \textit{Proof:}\quad Denoting $\textbf{w}:=\textbf{u}-\textbf{u}^1_\theta$ and $\textbf{H}:=\textbf{B}-\textbf{B}^1_\theta$,  we can obtain the following error equations:
 \begin{align}\label{eq23}
\frac{d\textbf{w}}{dt}+\nu A_f\textbf{w}+B[\textbf{u},\textbf{u}]-B[\textbf{u}^1_\theta,\textbf{u}^1_\theta]+SC_f[\textbf{B},\textbf{B}]-SC_f[\textbf{B}^1_\theta,\textbf{B}^1_\theta]=\Lambda,\\\label{eq24}
\frac{d\textbf{H}}{dt}+\mu A_B\textbf{H}-C_B[\textbf{u},\textbf{B}]+C_B[\textbf{u}^1_\theta,\textbf{B}^1_\theta]=\Pi.
 \end{align}
Since
 \begin{align}\label{eq25}
B[\textbf{u},\textbf{u}]-B[\textbf{u}^1_\theta,\textbf{u}^1_\theta]=B[\textbf{u},\textbf{w}]+B[\textbf{w},\textbf{u}]-B[\textbf{w},\textbf{w}],\\\label{eq26}
SC_f[\textbf{B},\textbf{B}]-SC_f[\textbf{B}^1_\theta,\textbf{B}^1_\theta]=SC_f[\textbf{B},\textbf{H}]+SC_f[\textbf{H},\textbf{B}]-SC_f[\textbf{H},\textbf{H}],\\\label{eq27}
C_B[\textbf{u},\textbf{B}]-C_B[\textbf{u}^1_\theta,\textbf{B}^1_\theta]=C_B[\textbf{u},\textbf{H}]+C_B[\textbf{w},\textbf{B}]-C_B[\textbf{w},\textbf{H}].
 \end{align}
 Taking inner product of \eqref{eq23}-\eqref{eq24} with $(\textbf{w},\textbf{H})$, one finds that
  \begin{align}\label{eq28}
&\frac{1}{2}\frac{d}{dt}\|\textbf{w}\|^2_{L^2(\Omega)}+\frac{S}{2}\frac{d}{dt}\|\textbf{H}\|^2_{L^2(\Omega)}+\nu\|\nabla\textbf{w}\|^2_{L^2(\Omega)}+\mu S\|\curl\textbf{H}\|^2_{L^2(\Omega)}\\&\nonumber\quad
+b(\textbf{w},\textbf{u},\textbf{w})+Sc_{\widetilde{B}}(\textbf{H},\textbf{B},\textbf{w})-Sc_{\widehat{B}}(\textbf{w},\textbf{B},\textbf{H})=-(\Lambda,\textbf{w})-S(\Pi,\textbf{H}).
 \end{align}
 Then it follows that
   \begin{align}\label{eq30}
&\frac{1}{2}\frac{d}{dt}\|\textbf{w}\|^2_{L^2(\Omega)}+\frac{S}{2}\frac{d}{dt}\|\textbf{H}\|^2_{L^2(\Omega)}+\nu\|\nabla\textbf{w}\|^2_{L^2(\Omega)}
+S\mu\|\curl\textbf{H}\|^2_{L^2(\Omega)}\\&\nonumber
\leq|(\Lambda,\textbf{w})|+S|(\Pi,\textbf{H})|+|b(\textbf{w},\textbf{u},\textbf{w})|+|Sc_{\widetilde{B}}(\textbf{H},\textbf{B},\textbf{w})|+S|c_{\widehat{B}}(\textbf{w},\textbf{B},\textbf{H})|.
 \end{align}
Thanks to the H\"{o}lder, Young and Embedded inequalities, we bound as follows:
    \begin{align*}
|(\Lambda,\textbf{w})|&\leq C\|\Lambda\|^2_{H^{-1}}+\frac{\nu}{8}\|\nabla\textbf{w}\|^2_{L^2(\Omega)},\\
|S(\Pi,\textbf{H})|&\leq C\|\Pi\|^2_{\cW^{-1}}+\frac{\mu S}{8}\|\curl\textbf{H}\|^2_{L^2(\Omega)},\\
|b(\textbf{w},\textbf{u},\textbf{w})|&\leq C\|\textbf{w}\|^2_{L^2(\Omega)}\|\nabla\textbf{u}\|^2_{L^2(\Omega)}+\frac{\nu}{8}\|\nabla\textbf{w}\|^2_{L^2(\Omega)},\\
|Sc_{\widetilde{B}}(\textbf{H},\textbf{B},\textbf{w})|&\leq C\|\textbf{H}\|^2_{L^2(\Omega)}\|\curl\textbf{B}\|^2_{L^2(\Omega)}+\frac{\nu}{8}\|\nabla\textbf{w}\|^2_{L^2(\Omega)},\\
|Sc_{\widehat{B}}(\textbf{w},\textbf{B},\textbf{H})|&\leq C\|\textbf{w}\|^2_{L^2(\Omega)}\|\curl\textbf{B}\|^2_{L^2(\Omega)}+\frac{\mu S}{8}\|\curl\textbf{H}\|^2_{L^2(\Omega)}.
 \end{align*}
Combining the above inequalities with \eqref{eq30},  we derive
    \begin{align}\label{eq31}
\frac{d}{dt}\bigl(\|\textbf{w}\|^2_{L^2(\Omega)}+\|\textbf{H}\|^2_{L^2(\Omega)}\bigr)
&-C\bigl(\|\nabla\textbf{u}\|^2_{L^2(\Omega)}+\|\curl\textbf{B}\|^2_{L^2(\Omega)}\bigr)\\\nonumber&
\times\bigl(\|\textbf{w}\|^2_{L^2(\Omega)}+\|\textbf{H}\|^2_{L^2(\Omega)}\bigr)
 \leq C\bigl(\|\Lambda\|^2_{H^{-1}}+\|\Pi\|^2_{\cW^{-1}}\bigr).
 \end{align}
 Applying the Gr\"{o}nwall inequality, one finds that
     \begin{align}\label{eq32}
&\|\textbf{w}\|^2_{L^2(\Omega)}+\|\textbf{H}\|^2_{L^2(\Omega)}\\&\nonumber
\leq \exp\bigl[\int_0^tC\bigl(\|\nabla\textbf{u}\|^2_{L^2(\Omega)}+\|\curl\textbf{B}\|^2_{L^2(\Omega)}\bigr)ds\bigr]\bigl(\|\textbf{w}(0)\|^2_{L^2(\Omega)}+\|\textbf{H}(0)\|^2_{L^2(\Omega)}\bigr)\\&\nonumber\quad
+\exp\bigl[\int_0^tC\bigl(\|\nabla\textbf{u}\|^2_{L^2(\Omega)}+\|\curl\textbf{B}\|^2_{L^2(\Omega)}\bigr)ds\bigr]\int_0^t\exp\bigl[-\int_0^sC\bigl(\|\nabla\textbf{u}\|^2_{L^2(\Omega)}\\&\nonumber\quad
+\|\curl\textbf{B}\|^2_{L^2(\Omega)}\bigr)d\tau\bigr]\bigl(\|\Lambda\|^2_{H^{-1}}+\|\Pi\|^2_{\cW^{-1}}\bigr)ds.
 \end{align}
 Moreover, we have
 \begin{align*}
 \exp\bigl[\int_0^tC\bigl(\|\nabla\textbf{u}\|^2_{L^2(\Omega)}+\|\curl\textbf{B}\|^2_{L^2(\Omega)}\bigr)ds\bigr]&\leq C\\\nonumber
\exp\bigl[-\int_0^sC\bigl(\|\nabla\textbf{u}\|^2_{L^2(\Omega)} +\|\curl\textbf{B}\|^2_{L^2(\Omega)}\bigr)d\tau\bigr] &\leq 1,\\\nonumber
\|\textbf{w}(0)\|^2_{L^2(\Omega)}+\|\textbf{H}(0)\|^2_{L^2(\Omega)}&\leq C\epsilon^2.
  \end{align*}
 Then we obtain
  \begin{align}\label{eq33}
\sup_{t\in [0,T]}\bigl(\|\textbf{u}(t)-\textbf{u}^1_{\theta}(t)\|^2_{L^2(\Omega)}+\|\textbf{B}(t)-\textbf{B}^1_{\theta}(t)\|^2_{L^2(\Omega)}\bigr)\leq C\epsilon.
 \end{align}
  The proof is completed.
 $$\eqno\Box$$

   \begin{lemma}%\cite{Biswas2022}
    Assume that
 \begin{align}\label{eq34}
&\|\textbf{u}^1_\theta|_{\partial\Omega}\|^4_{L^4([0,T], H^{\frac{1}{2}}(\partial\Omega))}+\|\textbf{B}^1_\theta\cdot\textbf{n}|_{\partial\Omega}\|^4_{L^4([0,T], H^{\frac{1}{2}}(\partial\Omega))}\\&\nonumber
\|\nabla\cdot\textbf{u}_\theta\|^4_{L^4([0,T], L^2(\Omega))}+\|\nabla\cdot\textbf{B}_\theta\|^4_{L^4([0,T], L^2(\Omega))}\leq C\epsilon^2,
 \end{align}
and applying the Hodge decomposition \eqref{eq42}-\eqref{eq43}, we have
  \begin{align}\label{eq35}
\|\textbf{u}^2_\theta\|^4_{L^4([0,T], \cX)}+\|\textbf{B}^2_\theta\|^4_{L^4([0,T], \cW)}\leq C\epsilon^2.
 \end{align}
\end{lemma}
\noindent \textit{Proof:}\quad Using the similar lines in \cite{Biswas2022}, here we omit its proof.
 $$\eqno\Box$$

\begin{theorem} Assume that $(\textbf{u},\textbf{B},p)$ is a strong solution of \eqref{eq1.1}-\eqref{eq3} and  $(\textbf{u}_\theta,\textbf{B}_\theta,p_\theta)$
  such that
 \begin{align}\label{eq36}
&\|\textbf{u}_\theta|_{\partial\Omega}\|^4_{L^4([0,T], H^{\frac{1}{2}}(\partial\Omega))}+
\|\textbf{u}_{\theta,0}-\textbf{u}_0\|^2_{L^2(\Omega)}\\&\nonumber\quad
+\|\mathfrak{L}_f[\textbf{u}_\theta,\textbf{B}_\theta,p_\theta]\|^2_{L^2([0,T]\times\Omega)}
+\|\mathfrak{L}_B[\textbf{u}_\theta,\textbf{B}_\theta]\|^2_{L^2([0,T]\times\Omega)}\\&\nonumber\quad
+ \|\textbf{B}_{\theta,0}-\textbf{B}_0\|^2_{L^2(\Omega)}+\|\textbf{B}_\theta\cdot\textbf{n}|_{\partial\Omega}\|^4_{L^4([0,T], H^{\frac{1}{2}}(\partial\Omega))}\\&\nonumber\quad
+\|\nabla\cdot\textbf{u}_\theta\|^4_{L^4([0,T], L^2(\Omega))}+\|\nabla\cdot\textbf{B}_\theta\|^4_{L^4([0,T], L^2(\Omega))}\\&\nonumber\quad
+\|\textbf{u}_\theta\|^4_{L^4([0,T], H^1(\Omega))}+\|\textbf{B}_\theta\|^4_{L^4([0,T], H^1(\Omega))}\leq C\epsilon^2.
 \end{align}
Then we have
 \begin{align}\label{eq37}
\|\textbf{u}-\textbf{u}_{\theta}\|^4_{L^2(\Omega_T)}+\|\textbf{B}-\textbf{B}_{\theta}\|^4_{L^2(\Omega_T)}\leq C\epsilon^2.
 \end{align}
 \end{theorem}
\noindent \textit{Proof:}\quad Let $(\textbf{u}^1_\theta, \textbf{B}^1_\theta)$ satisfy the following equations:
 \begin{align}\label{eq44}
\frac{d\textbf{u}^1_\theta}{dt}&+\nu A_f\textbf{u}^1_\theta+B[\textbf{u}_\theta,\textbf{u}_\theta]-B[\textbf{u}^1_\theta,\textbf{u}^1_\theta]
\\&\nonumber+SC_f[\textbf{B},\textbf{B}]
-SC_f[\textbf{B}^1_\theta,\textbf{B}^1_\theta]=\textbf{f}+\mathbb{P}\widehat{\textbf{f}},\\\label{eq45}
\frac{d\textbf{B}^1_\theta}{dt}&+\mu A_B\textbf{B}^1_\theta-C_B[\textbf{u}_\theta,\textbf{B}_\theta]+C_B[\textbf{u}^1_\theta,\textbf{B}^1_\theta]=0,
 \end{align}
For nonlinear term,  adding and subtracting some terms, we can rewrite,
 \begin{align}\label{eq46}
B[\textbf{u}_\theta,\textbf{u}_\theta]-B[\textbf{u}^1_\theta,\textbf{u}^1_\theta]=B[\textbf{u}^2_\theta,\textbf{u}_\theta]+B[\textbf{u}^1_\theta,\textbf{u}^2_\theta]:=\Psi_1,\\\label{eq47}
C_f[\textbf{B}_\theta,\textbf{B}_\theta]-C_f[\textbf{B}^1_\theta,\textbf{B}^1_\theta]=C_f[\textbf{B}^2_\theta,\textbf{B}_\theta]+C_f[\textbf{B}^1_\theta,\textbf{B}^2_\theta]:=\Psi_2,\\\label{eq48}
C_B[\textbf{u}_\theta,\textbf{B}_\theta]-C_B[\textbf{u}^1_\theta,\textbf{B}^1_\theta]=C_B[\textbf{u}^2_\theta,\textbf{B}_\theta]+C_B[\textbf{u}^1_\theta,\textbf{B}^2_\theta]:=\Psi_3.
 \end{align}
We will estimate $\int_0^T\|\Psi_1\|^2_{H^{-1}}dt$, $\int_0^T\|\Psi_2\|^2_{H^{-1}}dt$ and $\int_0^T\|\Psi_3\|^2_{\cW^{-1}}dt$ as follows.
Note that
 \begin{align}\label{eq49}
\|\Psi_1\|_{H^{-1}}&=\sup_{\textbf{w}\in \cX, \|\nabla\textbf{w}\|_{L^2(\Omega)}\leq 1}\bigl[b(\textbf{u}^2_\theta,\textbf{u}_\theta,\textbf{w})+b(\textbf{u}^1_\theta,\textbf{u}^2_\theta,\textbf{w})\bigr],\\\label{eq50}
\|\Psi_2\|_{H^{-1}}&=\sup_{\textbf{w}\in \cX, \|\nabla\textbf{w}\|_{L^2(\Omega)}\leq 1}\bigl[c_{\widetilde{B}}(\textbf{B}^2_\theta,\textbf{B}_\theta,\textbf{w})+c_{\widetilde{B}}(\textbf{B}^1_\theta,\textbf{B}^2_\theta,\textbf{w})\bigr],\\\label{eq51}
\|\Psi_3\|_{\cW^{-1}}&=\sup_{\textbf{H}\in \cW, \|\curl\textbf{H}\|_{L^2(\Omega)}\leq 1}\bigl[c_{\widehat{B}}(\textbf{u}^2_\theta,\textbf{B}_\theta,\textbf{H})+c_{\widehat{B}}(\textbf{u}^1_\theta,\textbf{B}^2_\theta,\textbf{H})\bigr].
 \end{align}
 We bound the term \eqref{eq49}-\eqref{eq51} as follows:
\begin{align*}
b(\textbf{u}^2_\theta,\textbf{u}_\theta,\textbf{w})&\leq C\|\textbf{u}^2_\theta\|_{L^4(\Omega)}\|\nabla\textbf{u}_\theta\|_{L^2(\Omega)}\|\textbf{w}\|_{L^4(\Omega)}\\\nonumber&
\leq C\|\textbf{u}^2_\theta\|^{\frac{1}{2}}_{L^2(\Omega)}\|\nabla\textbf{u}^2_\theta\|^{\frac{1}{2}}_{L^2(\Omega)}\|\nabla\textbf{u}_\theta\|_{L^2(\Omega)}\|\nabla\textbf{w}\|_{L^2(\Omega)},\\\nonumber
b(\textbf{u}^1_\theta,\textbf{u}^2_\theta,\textbf{w})&\leq C\|\textbf{u}^1_\theta\|_{L^4(\Omega)}\|\nabla\textbf{u}^2_\theta\|_{L^2(\Omega)}\|\textbf{w}\|_{L^4(\Omega)}\\\nonumber&
\leq C\|\textbf{u}^1_\theta\|_{H^1(\Omega)}\|\nabla\textbf{u}^2_\theta\|_{L^2(\Omega)}\|\nabla\textbf{w}\|_{L^2(\Omega)},
\end{align*}
\begin{align*}
Sc_{\widetilde{B}}(\textbf{B}^2_\theta,\textbf{B}_\theta,\textbf{w})&\leq C\|\textbf{B}^2_\theta\|_{L^4(\Omega)}\|\nabla\textbf{B}_\theta\|_{L^2(\Omega)}\|\textbf{w}\|_{L^4(\Omega)}\\\nonumber&
\leq C\|\textbf{B}^2_\theta\|^{\frac{1}{2}}_{L^2(\Omega)}\|\nabla\textbf{B}^2_\theta\|^{\frac{1}{2}}_{L^2(\Omega)}\|\nabla\textbf{B}_\theta\|_{L^2(\Omega)}\|\nabla\textbf{w}\|_{L^2(\Omega)},\\\nonumber
Sc_{\widetilde{B}}(\textbf{B}^1_\theta,\textbf{B}^2_\theta,\textbf{w})&\leq C\|\textbf{B}^1_\theta\|_{L^4(\Omega)}\|\nabla\textbf{B}^2_\theta\|_{L^2(\Omega)}\|\textbf{w}\|_{L^4(\Omega)}\\\nonumber&
\leq C\|\textbf{B}^1_\theta\|_{H^1(\Omega)}\|\nabla\textbf{B}^2_\theta\|_{L^2(\Omega)}\|\nabla\textbf{w}\|_{L^2(\Omega)},
\end{align*}
and
\begin{align*}
c_{\widehat{B}}(\textbf{u}^2_\theta,\textbf{B}_\theta,\textbf{H})&\leq C\|\textbf{u}^2_\theta\|_{L^4(\Omega)}\|\textbf{B}_\theta\|_{L^4(\Omega)}\|\curl\textbf{H}\|_{L^2(\Omega)}\\\nonumber&
\leq C\|\textbf{u}^2_\theta\|^{\frac{1}{2}}_{L^2(\Omega)}\|\nabla\textbf{u}^2_\theta\|^{\frac{1}{2}}_{L^2(\Omega)}\|\textbf{B}_\theta\|_{H^1(\Omega)}\|\curl\textbf{H}\|_{L^2(\Omega)},\\\nonumber
c_{\widehat{B}}(\textbf{u}^1_\theta,\textbf{B}^2_\theta,\textbf{H})&\leq C\|\textbf{u}^1_\theta\|_{L^4(\Omega)}\|\nabla\textbf{B}^2_\theta\|_{L^2(\Omega)}\|\curl\textbf{H}\|_{L^2(\Omega)}\\\nonumber&
\leq C\|\nabla\textbf{u}^1_\theta\|_{L^2(\Omega)}\|\textbf{B}^2_\theta\|_{H^1(\Omega)}\|\curl\textbf{H}\|_{L^2(\Omega)}.
 \end{align*}
 Then we have
  \begin{align}\label{eq52}
&\int_{0}^T\|\Psi_1\|_{H^{-1}}^2dt\\\nonumber
&\leq C\int_{0}^T\|\nabla\textbf{u}_\theta\|^2_{L^2(\Omega)}\|\textbf{u}^2_\theta\|_{L^2(\Omega)}\|\nabla\textbf{u}^2_\theta\|_{L^2(\Omega)}dt\\&\nonumber\quad+
C\int_0^T\|\nabla\textbf{u}_\theta^2\|^2_{L^2(\Omega)}\|\nabla\textbf{u}^1_\theta\|^2_{L^2(\Omega)}dt\\&\nonumber
\leq C\Bigl(\int_0^T\|\nabla\textbf{u}_\theta\|^4_{L^2(\Omega)}dt\Bigr)^{\frac{1}{2}}\Bigl(\int_0^T\|\nabla\textbf{u}_\theta^2\|^4_{L^2(\Omega)}dt\Bigr)^{\frac{1}{4}}\Bigl(\int_0^T\|\textbf{u}_\theta^2\|^4_{L^2(\Omega)}dt\Bigr)^{\frac{1}{4}}\\&\nonumber\quad
+C\Bigl(\int_0^T\|\nabla\textbf{u}^2_\theta\|^4_{L^2(\Omega)}dt\Bigr)^{\frac{1}{2}}\Bigl(\int_0^T\|\nabla\textbf{u}_\theta^1\|^4_{H^1(\Omega)}dt\Bigr)^{\frac{1}{2}}\\&\nonumber
\leq C\epsilon,
\end{align}
  \begin{align}\label{eq53}
&\int_{0}^T\|\Psi_2\|_{H^{-1}}^2dt\\\nonumber
&\leq C\int_{0}^T\|\nabla\textbf{B}_\theta\|^2_{L^2(\Omega)}\|\textbf{B}^2_\theta\|_{L^2(\Omega)}\|\nabla\textbf{B}^2_\theta\|_{L^2(\Omega)}dt\\&\nonumber\quad+
C\int_0^T\|\nabla\textbf{B}_\theta^2\|^2_{L^2(\Omega)}\|\textbf{B}^1_\theta\|^2_{H^1(\Omega)}dt\\&\nonumber
\leq C\Bigl(\int_0^T\|\nabla\textbf{B}_\theta\|^4_{L^2(\Omega)}dt\Bigr)^{\frac{1}{2}}\Bigl(\int_0^T\|\nabla\textbf{B}_\theta^2\|^4_{L^2(\Omega)}dt\Bigr)^{\frac{1}{4}}\Bigl(\int_0^T\|\textbf{B}_\theta^2\|^4_{L^2(\Omega)}dt\Bigr)^{\frac{1}{4}}\\&\nonumber\quad
+C\Bigl(\int_0^T\|\nabla\textbf{B}^2_\theta\|^4_{L^2(\Omega)}dt\Bigr)^{\frac{1}{2}}\Bigl(\int_0^T\|\textbf{B}_\theta^1\|^4_{H^1(\Omega)}dt\Bigr)^{\frac{1}{2}}\\&\nonumber
\leq C\epsilon,
\end{align}
and
  \begin{align}\label{eq54}
&\int_{0}^T\|\Psi_3\|_{\cW^{-1}}^2dt\\\nonumber
&\leq C\int_{0}^T\|\nabla\textbf{u}^2_\theta\|_{L^2(\Omega)}\|\textbf{u}^2_\theta\|_{L^2(\Omega)}\|\textbf{B}_\theta\|^2_{H^1(\Omega)}dt\\&\nonumber\quad+
C\int_0^T\|\nabla\textbf{u}_\theta^1\|^2_{L^2(\Omega)}\|\textbf{B}^2_\theta\|^2_{H^1(\Omega)}dt\\&\nonumber
\leq C\Bigl(\int_0^T\|\nabla\textbf{u}^2_\theta\|^4_{L^2(\Omega)}dt\Bigr)^{\frac{1}{4}}\Bigl(\int_0^T\|\nabla\textbf{u}^2_\theta\|^4_{L^2(\Omega)}dt\Bigr)^{\frac{1}{4}}\Bigl(\int_0^T\|\textbf{B}_\theta\|^4_{H^1(\Omega)}dt\Bigr)^{\frac{1}{2}}\\&\nonumber\quad
+C\Bigl(\int_0^T\|\nabla\textbf{u}^1_\theta\|^4_{L^2(\Omega)}dt\Bigr)^{\frac{1}{2}}\Bigl(\int_0^T\|\textbf{B}_\theta^2\|^4_{H^1(\Omega)}dt\Bigr)^{\frac{1}{2}}
\\&\nonumber
\leq C\epsilon.
\end{align}
Moreover, since
  \begin{align}\label{eq55}
\|\mathbb{P}\Delta \textbf{u}_\theta^2\|_{H^{-1}}&\leq C\|\nabla\textbf{u}_\theta^2\|_{L^2(\Omega)}\|\textbf{w}\|_{H^1(\Omega)},\\\nonumber
\|\mathbb{Q}\curl\curl\textbf{B}_\theta^2\|_{\cW^{-1}}&\leq C\|\curl\textbf{B}_\theta^2\|_{L^2(\Omega)}\|\curl\textbf{H}\|_{L^2(\Omega)}.
\end{align}
Thus, it follows that
  \begin{align}\label{eq56}
\int_0^T\|\mathbb{P}\Delta \textbf{u}_\theta^2\|^2_{H^{-1}}dt&\leq C\int_0^T\|\nabla\textbf{u}_\theta^2\|^2_{L^2(\Omega)}\|\textbf{w}\|^2_{H^1(\Omega)}dt\\&\nonumber
\leq  C\Bigl(\int_0^T\|\nabla\textbf{u}_\theta^2\|^4_{L^2(\Omega)}dt\Bigr)^{\frac{1}{2}}\Bigl(\int_0^T\|\textbf{w}\|^4_{H^1(\Omega)}dt\Bigr)^{\frac{1}{2}}\\&\nonumber
\leq C\epsilon,\\\label{eq57}
\int_0^T\|\mathbb{Q}\curl\curl\textbf{B}_\theta^2\|^2_{\cW^{-1}}dt&\leq C\int_0^T\|\curl\textbf{B}_\theta^2\|^2_{L^2(\Omega)}\|\curl\textbf{H}\|^2_{L^2(\Omega)}dt\\&\nonumber
\leq  C\Bigl(\int_0^T\|\curl\textbf{B}_\theta^2\|^4_{L^2(\Omega)}dt\Bigr)^{\frac{1}{2}}\Bigl(\int_0^T\|\curl\textbf{H}\|^4_{L^2(\Omega)}dt\Bigr)^{\frac{1}{2}}\\&\nonumber
\leq C\epsilon.
\end{align}
Denoting $\phi:=\mathbb{P}\widehat{\textbf{f}}+\mathbb{P}(\Delta\textbf{u}_\theta^2)-\Psi_1-\Psi_2$ and $\psi:=\mathbb{Q}(\curl\curl\textbf{B}_\theta^2)-\Psi_3$, we have
 \begin{align}\label{eq58}
\int_0^T\|\phi\|^2_{H^{-1}}dt&\leq C\epsilon,\\\label{eq59}
\int_0^T\|\psi\|^2_{\cW^{-1}}dt&\leq C\epsilon.
 \end{align}
 Then we have
 \begin{align}\label{eq60}
\frac{d\textbf{u}^1_\theta}{dt}+\nu A_f\textbf{u}^1_\theta+B[\textbf{u}^1_\theta,\textbf{u}^1_\theta]+SC_f[\textbf{B}^1_\theta,\textbf{B}^1_\theta]&=\mathbb{P}\widehat{\textbf{f}}+\phi,\\\label{eq61}
\frac{d\textbf{B}^1_\theta}{dt}+\mu A_B\textbf{B}^1_\theta-C_B[\textbf{u}^1_\theta,\textbf{B}^1_\theta]&=\psi,
 \end{align}
 Applying \textbf{Theorem 3}, we have
  \begin{align}\label{eq62}
\sup_{t \in [0,T]}\|\textbf{u}(t)-\textbf{u}^1_\theta(t)\|_{L^2(\Omega)}&\leq C\epsilon,\\\label{eq63}
\sup_{t \in [0,T]}\|\textbf{B}(t)-\textbf{B}^1_\theta(t)\|_{L^2(\Omega)}&\leq C\epsilon.
 \end{align}
 Moreover, since
 \begin{align}\label{eq64}
\Bigl(\int_0^T\|\textbf{u}_\theta^2\|^4_{L^2}dt\Bigr)^{\frac{1}{4}}+\Bigl(\int_0^T\|\textbf{B}_\theta^2\|^4_{L^2}dt\Bigr)^{\frac{1}{4}}\leq C\sqrt{\epsilon}.
 \end{align}
  Then we have
   \begin{align}\label{eq65}
\int_0^T\|\textbf{u}-\textbf{u}_\theta\|^4_{L^2}dt\leq \int_0^T\|\textbf{u}-\textbf{u}^1_\theta\|^4_{L^2}dt+\int_0^T\|\textbf{u}^2_\theta\|^4_{L^2}dt\leq C\epsilon^2,\\\label{eq66}
\int_0^T\|\textbf{B}-\textbf{B}_\theta\|^4_{L^2}dt\leq \int_0^T\|\textbf{B}-\textbf{B}^1_\theta\|^4_{L^2}dt+\int_0^T\|\textbf{B}^2_\theta\|^4_{L^2}dt\leq C\epsilon^2.
  \end{align}
 The proof is completed.
 $$\eqno\Box$$

\begin{theorem} Given $\epsilon>0$, we can find $(\textbf{u}_\theta,\textbf{B}_\theta,p_\theta)\in \mathfrak{F}_N$, such that
 \begin{align}\label{eq66}
&\|\textbf{u}_\theta|_{\partial\Omega}\|^4_{L^4([0,T], H^{\frac{1}{2}}(\partial\Omega))}+
\|\textbf{u}_{\theta,0}-\textbf{u}_0\|^2_{L^2(\Omega)}\\&\quad\nonumber
+\|\mathfrak{L}_f[\textbf{u}_\theta,\textbf{B}_\theta,p_\theta]\|^2_{L^2(\Omega_T)}
+\|\mathfrak{L}_B[\textbf{u}_\theta,\textbf{B}_\theta]\|^2_{L^2(\Omega_T)}\\&\quad\nonumber
+ \|\textbf{B}_{\theta,0}-\textbf{B}_0\|^2_{L^2(\Omega)}+\|\textbf{B}_\theta\cdot\textbf{n}|_{\partial\Omega}\|^4_{L^4([0,T], H^{\frac{1}{2}}(\partial\Omega))}\\&\quad\nonumber
+\|\nabla\cdot\textbf{u}_\theta\|^4_{L^4([0,T], L^2(\Omega))}+\|\nabla\cdot\textbf{B}_\theta\|^4_{L^4([0,T], L^2(\Omega))}\\&\quad\nonumber
+\|\textbf{u}_\theta\|^4_{L^4([0,T], H^1(\Omega))}+\|\textbf{B}_\theta\|^4_{L^4([0,T], H^1(\Omega))}\leq C\epsilon^2.
  \end{align}
Moveover, our scheme is approximately stable.
 \end{theorem}
\noindent \textit{Proof:}\quad From \textbf{Lemma 1}, given $\epsilon>0$, let $(\textbf{u},\textbf{B},p)$ be a strong solution of problem \eqref{eq1.1}-\eqref{eq3} and
$(\textbf{u}_\theta,\textbf{B}_\theta,p_\theta)\in \mathfrak{F}_N$ satisfying
\begin{subequations}\label{eq68.1}
\begin{alignat}{2} \label{eq68.1a}
&\sup_{t \in [0,T]}\|\textbf{u}(t)-\textbf{u}_\theta(t)\|_{L^2(\Omega)}\leq C\epsilon,\\\label{eq68.1b}
&\sup_{t \in [0,T]}\|\textbf{B}(t)-\textbf{B}_\theta(t)\|_{L^2(\Omega)}\leq C\epsilon,\\\label{eq68.1c}
&\|\textbf{u}-\textbf{u}_\theta\|_{H^{1,2}(\Omega_T)}\leq C\epsilon,\\\label{eq68.1d}
&\|\textbf{B}-\textbf{B}_\theta\|_{H^{1,2}(\Omega_T)}\leq C\epsilon,\\\label{eq68.1e}
&\|\textbf{u}-\textbf{u}_\theta\|_{L^4([0,T]\times W^{1,4}(\Omega))}\leq C\epsilon,\\\label{eq68.1f}
&\|\textbf{B}-\textbf{B}_\theta\|_{L^4([0,T]\times W^{1,4}(\Omega))}\leq C\epsilon,\\\label{eq68.1g}
&\|p-p_\theta\|_{L^2([0,T]\times H^1(\Omega))}\leq C\epsilon.
\end{alignat}
\end{subequations}
By  \eqref{eq68.1a} and \eqref{eq68.1b}, we have
\begin{align}\label{eq69}
\|\textbf{u}(0)-\textbf{u}_\theta(0)\|^2_{L^2(\Omega)}+\|\textbf{B}(0)-\textbf{B}_\theta(0)\|^2_{L^2(\Omega)}\leq C\epsilon^2.
 \end{align}
Using  \eqref{eq68.1e} and \eqref{eq68.1f}, we obtain
\begin{align}\label{eq70}
&\|\nabla\cdot\textbf{u}\|^4_{L^4([0,T],L^2(\Omega))}+\|\nabla\cdot\textbf{B}\|^4_{L^4([0,T],L^2(\Omega))}\\&\nonumber
=\|\nabla\cdot\textbf{u}-\nabla\cdot\textbf{u}_\theta\|^4_{L^4([0,T],L^2(\Omega))}+\|\nabla\cdot\textbf{B}-\nabla\cdot\textbf{B}_\theta\|^4_{L^4([0,T],L^2(\Omega))}\\&\nonumber
\leq\|\textbf{u}-\textbf{u}_\theta\|^4_{L^4([0,T],H^1(\Omega))}+\|\textbf{B}-\textbf{B}_\theta\|^4_{L^4([0,T],H^1(\Omega))}\\&\nonumber
\leq C\epsilon^4.
 \end{align}
Taking $\gamma>0$ small enough, one finds that
 \begin{align}\label{eq71}
&\quad\gamma\bigl(\|\textbf{u}_\theta\|^4_{L^4([0,T],H^1(\Omega))}+\|\textbf{B}_\theta\|^4_{L^4([0,T],H^1(\Omega))}\bigr)\\&\nonumber
\leq \gamma C\bigl(\|\textbf{u}-\textbf{u}_\theta\|^4_{L^4([0,T],H^1(\Omega))}+\|\textbf{B}-\textbf{B}_\theta\|^4_{L^4([0,T],H^1(\Omega))}\bigr)\\&\nonumber\quad+
\gamma C\bigl(\|\textbf{u}\|^4_{L^4([0,T],H^1(\Omega))}+\|\textbf{B}\|^4_{L^4([0,T],H^1(\Omega))}\bigr)\\&\nonumber
\leq C\epsilon^4.
 \end{align}
 Consider
 \begin{align}\label{eq72}
&\|\mathfrak{L}_f[\textbf{u}_\theta,\textbf{B}_\theta, p_\theta]-\mathfrak{L}_f[\textbf{u},\textbf{B}, p]\|^2_{L^2(\Omega_T)}
+\|\mathfrak{L}_B[\textbf{u}_\theta,\textbf{B}_\theta]-\mathfrak{L}_B[\textbf{u},\textbf{B}]\|^2_{L^2(\Omega_T)}\\\nonumber&
\leq C\|\partial_t\textbf{u}-\partial_t\textbf{u}_\theta\|^2_{L^2(\Omega_T)}+C\|\partial_t\textbf{B}-\partial_t\textbf{B}_\theta\|^2_{L^2(\Omega_T)}\\\nonumber&\quad+
C\|\Delta\textbf{u}-\Delta\textbf{u}_\theta\|^2_{L^2(\Omega_T)}+C\|\curl\curl(\textbf{B}-\textbf{B}_\theta)\|^2_{L^2(\Omega_T)}\\\nonumber&\quad+
C\|\nabla p-\nabla p_\theta\|^2_{L^2(\Omega_T)}+C\|B[\textbf{u},\textbf{u}]-B[\textbf{u}_\theta,\textbf{u}_\theta]\|^2_{L^2(\Omega_T)}\\\nonumber&\quad+
C\|C_f[\textbf{B},\textbf{B}]-C_f[\textbf{B}_\theta,\textbf{B}_\theta]\|^2_{L^2(\Omega_T)}\\\nonumber&\quad
+C\|C_B[\textbf{u},\textbf{B}]-C_B[\textbf{u}_\theta,\textbf{B}_\theta]\|^2_{L^2(\Omega_T)}.
 \end{align}
 Using \eqref{eq76}-\eqref{eq78}, we arrive at
   \begin{align}\label{eq82}
 & \|B[\textbf{u},\textbf{u}]-B[\textbf{u}_\theta,\textbf{u}_\theta]\|^2_{L^2(\Omega_T)} \leq C\epsilon^2,\\\label{eq83}
&\|C_f[\textbf{B},\textbf{B}]-C_f[\textbf{B}_\theta,\textbf{B}_\theta]\|^2_{L^2(\Omega_T)}\leq C\epsilon^2,\\\label{eq84}
&\|C_B[\textbf{u},\textbf{B}]-C_B[\textbf{u}_\theta,\textbf{B}_\theta]\|^2_{L^2(\Omega_T)}\leq C\epsilon^2.
 \end{align}
Thanks to \eqref{eq73}-\eqref{eq81} and \eqref{eq82}-\eqref{eq84}, we have
  \begin{align}\label{eq85}
&\|\mathfrak{L}_f[\textbf{u}_\theta,\textbf{B}_\theta, p_\theta]-\mathfrak{L}_f[\textbf{u},\textbf{B}, p]\|^2_{L^2(\Omega_T)}
+\|\mathfrak{L}_B[\textbf{u}_\theta,\textbf{B}_\theta]-\mathfrak{L}_B[\textbf{u},\textbf{B}]\|^2_{L^2(\Omega_T)}\leq C\epsilon^2.
 \end{align}
Using \eqref{eq68.1e}-\eqref{eq68.1f} and the Triangle inequality, we derive
  \begin{align}\label{eq86}
&\|\textbf{u}|_{\partial\Omega}\|^4_{L^4([0,T], H^{\frac{1}{2}}\partial\Omega))}+\|\textbf{B}\cdot\textbf{n}|_{\partial\Omega}\|^4_{L^4([0,T], H^{\frac{1}{2}}\partial\Omega))}\\\nonumber&
=\|\textbf{u}|_{\partial\Omega}-\textbf{u}_\theta|_{\partial\Omega}\|^4_{L^4([0,T], H^{\frac{1}{2}}\partial\Omega))}+\|\textbf{B}\cdot\textbf{n}|_{\partial\Omega}-\textbf{B}_\theta\cdot\textbf{n}|_{\partial\Omega}\|^4_{L^4([0,T], H^{\frac{1}{2}}\partial\Omega))}\\\nonumber&\leq
\|\textbf{u}-\textbf{u}_\theta\|^4_{L^4([0,T], H^{1}(\Omega))}+\|\textbf{B}-\textbf{B}_\theta\|^4_{L^4([0,T], H^{1}(\Omega))}\\\nonumber&\leq
C\epsilon^4.
 \end{align}
 We obtain the desired result. The proof is finished.
 $$\eqno\Box$$

 \begin{theorem} Assume that $(\textbf{u}_{\theta,1},\textbf{B}_{\theta,1},p_{\theta,1})\in\mathfrak{F}_N$ is the approximate solution of
\begin{subequations}\label{eq87.1}
\begin{alignat}{2} \label{eq87.1a} \partial_t \textbf{u}-\nu\Delta\textbf{u}
+(\textbf{u}\cdot \nabla)\textbf{u}\\\nonumber
+S\textbf{B}\times\mbox{curl}\textbf{B}+\nabla p&=\textbf{f}_1,   \\
\partial_t \textbf{B}+\mu\mbox{curl}\mbox{curl} \textbf{B}-\mbox{curl}(\textbf{u}\times \textbf{B})&=0, \\
\mbox{div} \textbf{u}&=0, \\
\mbox{div} \textbf{B}&=0,\\
\textbf{u}{|_{\partial\Omega}}&=\textbf{0},  \\
 \textbf{B}\cdot \textbf{n}{|_{\partial\Omega}}&=0, \\
\mbox{curl \textbf{B}}\times \textbf{n}{|_{\partial\Omega}}&=0, \\
\textbf{u}(\textbf{x},0)&=\textbf{u}_{0,1}(\textbf{x}),\\
\textbf{B}(\textbf{x},0)&=\textbf{B}_{0,1}(\textbf{x}).
\end{alignat}
\end{subequations}
 Assume that $(\textbf{u}_{\theta,2},\textbf{B}_{\theta,2},p_{\theta,2})\in\mathfrak{F}_N$ is the approximate solution of
\begin{subequations}\label{eq88.1}
\begin{alignat}{2} \label{eq88.1a} \partial_t \textbf{u}-\nu\Delta\textbf{u}
+(\textbf{u}\cdot \nabla)\textbf{u}\\\nonumber
+S\textbf{B}\times\mbox{curl}\textbf{B}+\nabla p&=\textbf{f}_2,   \\
\partial_t \textbf{B}+\mu\mbox{curl}\mbox{curl} \textbf{B}-\mbox{curl}(\textbf{u}\times \textbf{B})&=0, \\
\mbox{div} \textbf{u}&=0, \\
\mbox{div} \textbf{B}&=0,\\
\textbf{u}{|_{\partial\Omega}}&=\textbf{0},  \\
 \textbf{B}\cdot \textbf{n}{|_{\partial\Omega}}&=0, \\
\mbox{curl \textbf{B}}\times \textbf{n}{|_{\partial\Omega}}&=0, \\
\textbf{u}(\textbf{x},0)&=\textbf{u}_{0,2}(\textbf{x}),\\
\textbf{B}(\textbf{x},0)&=\textbf{B}_{0,2}(\textbf{x}).
\end{alignat}
\end{subequations}
Then we have
  \begin{align*}
&\quad\|\textbf{u}_{\theta,1}-\textbf{u}_{\theta,2}\|_{L^4([0,T], L^2(\Omega))}+\|\textbf{B}_{\theta,1}-\textbf{B}_{\theta,2}\|_{L^4([0,T], L^2(\Omega))}\\\nonumber&\leq
C\epsilon^{\frac{1}{2}}+C\|\textbf{u}_{0,1}-\textbf{u}_{0,2}\|_{L^2(\Omega)}+C\|\textbf{B}_{0,1}-\textbf{B}_{0,2}\|_{L^2(\Omega)}+C\|\textbf{f}_{1}-\textbf{f}_{2}\|_{L^2(0,T,L^2(\Omega))}.
 \end{align*}
Moveover, our scheme is approximately stable.
 \end{theorem}
 \noindent \textit{Proof:}\quad
Using the Triangle inequality, we obtain
\begin{align}\label{eq89}
&\quad\|\textbf{u}_{\theta,1}-\textbf{u}_{\theta,2}\|_{L^4([0,T], L^2(\Omega))}+\|\textbf{B}_{\theta,1}-\textbf{B}_{\theta,2}\|_{L^4([0,T], L^2(\Omega))}\\\nonumber&
 \leq \|\textbf{u}_{\theta,1}-\textbf{u}_1\|_{L^4([0,T], L^2(\Omega))}+\|\textbf{u}_{\theta,2}-\textbf{u}_{2}\|_{L^4([0,T], L^2(\Omega))}+\|\textbf{u}_{1}-\textbf{u}_{2}\|_{L^4([0,T], L^2(\Omega))}\\\nonumber&\quad
 +\|\textbf{B}_{\theta,1}-\textbf{B}_1\|_{L^4([0,T], L^2(\Omega))}+\|\textbf{B}_{\theta,2}-\textbf{B}_{2}\|_{L^4([0,T], L^2(\Omega))}+\|\textbf{B}_{1}-\textbf{B}_{2}\|_{L^4([0,T], L^2(\Omega))}\\\nonumber&\quad
 \end{align}
 Thanks to
 \begin{align}\label{eq90}
&\quad\|\textbf{u}_{1}-\textbf{u}_{\theta,1}\|_{L^4([0,T], L^2(\Omega))}+\|\textbf{B}_{1}-\textbf{B}_{\theta,1}\|_{L^4([0,T], L^2(\Omega))}
\leq C\epsilon^{\frac{1}{2}}.
 \end{align}
  and
  \begin{align}\label{eq91}
&\quad\|\textbf{u}_{2}-\textbf{u}_{\theta,2}\|_{L^4([0,T], L^2(\Omega))}+\|\textbf{B}_{2}-\textbf{B}_{\theta,2}\|_{L^4([0,T], L^2(\Omega))}
\leq C\epsilon^{\frac{1}{2}}.
 \end{align}
 And using the stability of MHD, we have
   \begin{align}\label{eq92}
&\quad\|\textbf{u}_{1}-\textbf{u}_{2}\|_{L^4([0,T], L^2(\Omega))}+\|\textbf{B}_{1}-\textbf{B}_{2}\|_{L^4([0,T], L^2(\Omega))}\\\nonumber&\leq
C\|\textbf{u}_{0,1}-\textbf{u}_{0,2}\|_{L^2(\Omega)}+C\|\textbf{B}_{0,1}-\textbf{B}_{0,2}\|_{L^2(\Omega)}+C\|\textbf{f}_{1}-\textbf{f}_{2}\|_{L^2(0,T,L^2(\Omega))}.
 \end{align}
 Combining \eqref{eq90}-\eqref{eq92} with \eqref{eq89}, the desired result is obtained. The proof is finished.
 $$\eqno\Box$$

 %\textbf{Remark 2:} (1) Some proofs of Theorems 2-4 extend from some techniques proposed
%in reference \cite{Tabata2005} to magnetohydrodynamics problem with temperature dependent parameters \eqref{eq1.1}-\eqref{eq3}.
%
%(2) From the proof of Theorem 4, we derive the error estimate under $\Delta t\leq Ch$ condition for
% the fluid pressure for Algorithm 3.1.
%
%(3) We emphasize that the error estimate for the fluid pressure is not optimal from numerical results in section 5.

 %%%%%%%%%%%%%%%%%%%%%%%%%%%%%%%%%%%%%%%%%%%%%%%%%%%%%%%%%%%%%%%%%%%%%%%%%%%%%%%%%%%%%%%%%%%%%%%%%%%%%%%%%%%%%%%%%%%%%

%
%\section*{\label{Sec6}Conclusions}
%
%In this paper, we have analyzed a fully discrete  scheme
%for computing Cahn-Hilliard-Magneto-hydrodynamics system.
% The scheme is based on using conforming finite element method
%in space and Euler semi-implicit discretization with convex splitting in time. We have prove our scheme is unconditionally energy stable and
%obtain optimal error estimates for the concentration field, the chemical potential, the velocity field, the magnetic field
%and the pressure. Numerical tests are shown to confirm the theoretical rates of the our scheme.


\end{document}